\documentclass[11pt]{article}

\usepackage[margin=1in]{geometry}
\usepackage[T1]{fontenc}
\usepackage{lmodern}
\usepackage{microtype}
\usepackage{amsmath,amssymb,amsthm,mathtools}
\usepackage{enumitem}
\usepackage{tikz}
\usetikzlibrary{calc}
\usetikzlibrary{arrows.meta,decorations.markings}

\usepackage{xcolor}
\usepackage{hyperref}
\usepackage{cleveref}
\usepackage{caption}
\usepackage{float}

\hypersetup{hidelinks}
\newtheorem{theorem}{Theorem}[section]
\newtheorem{lemma}[theorem]{Lemma}
\newtheorem{claim}[theorem]{Claim}
\newtheorem{fact}[theorem]{Fact}
\newtheorem{corollary}[theorem]{Corollary}
\theoremstyle{definition}
\newtheorem{definition}[theorem]{Definition}
\newtheorem{construction}[theorem]{Construction}
\theoremstyle{remark}
\newtheorem*{remark}{Remark}

\newcommand{\rinf}{\rho_{\infty}}
\newcommand{\linf}{L_{\infty}}
\newcommand{\eps}{\varepsilon}
\newcommand{\cP}{\mathcal{P}}
\newcommand{\RR}{\mathbb{R}}
\newcommand{\ZZ}{\mathbb{Z}}
\newcommand{\PL}{P_{\mathrm L}}
\newcommand{\PR}{P_{\mathrm R}}
\newcommand{\PB}{P_{\mathrm B}}
\newcommand{\PU}{P_{\mathrm U}}
\newcommand{\PH}{P_{\mathrm H}}
\newcommand{\PV}{P_{\mathrm V}}

\makeatletter
\newcommand{\inlinetag}[3]{%
 \begingroup
 \phantomsection
 \def\@currentlabel{#2}%
 \label{#3}%
 \endgroup
 \(#1\)\nobreak\hspace{0.25em}\textup{(#2)}%
}
\makeatother

\title{A Complete Dichotomy for $L_\infty$ Permutation Removal}
\author{
Viet Hung Le\thanks{Department of Mathematics, Duke University, Durham, NC 27710, USA. Email: \texttt{hung.le@duke.edu}.}
\and
Ruilin Shi\thanks{Department of Mathematics, Duke University, Durham, NC 27710, USA. Email: \texttt{ruilin.shi@duke.edu}.}
\and
Fan Wei\thanks{Department of Mathematics, Duke University, Durham, NC 27710, USA. Supported by NSF grant DMS-2401414. Email: \texttt{fan.wei@duke.edu}.}
}
\date{}

\begin{document}
\maketitle

\begin{abstract}
We study the removal problem for permutation patterns under the
normalized $L_\infty$ earth mover's distance. We prove a
complete dichotomy. Let
\[
\mathcal E
 =\{(12),(21),(132),(213),(231),(312)\}.
\]
If $\sigma\in\mathcal E$, then every permutation of length $n$ that
is $\varepsilon$-far from being $\sigma$-free contains at least
$\varepsilon(n-1)$ copies of $\sigma$; for
$0<\varepsilon\le 1/16$, this is sharp up to a constant factor for
all sufficiently large $n$. If $\sigma\in S_k\setminus\mathcal E$,
then for every sufficiently large $n$ there is a permutation of
length $n$ that is $1/(8k)$-far from being $\sigma$-free and contains
exactly one copy of $\sigma$.

Thus the usual positive-density removal conclusion fails for every
nontrivial permutation pattern, but the sparse behavior is completely
characterized: the minimum number of copies is linear for exactly six
patterns and is one for all others. Consequently, no canonical
one-sided tester that rejects only after finding a forbidden copy can
have query complexity independent of the permutation length. 

The results and proofs in this paper have also been formally
verified in Lean~4; see the ``Formal verification'' section
following the main text for details.
\end{abstract}

\section{Introduction}

The graph removal lemma is one of the main tools in extremal combinatorics. To state it, let $G$ and $H$ be graphs on the same vertex set $[n]$, and define their normalized edit distance by
$d_1(G,H)
 =\frac{|E(G)\mathbin{\triangle}E(H)|}{\binom{n}{2}}$, the fraction of vertex pairs whose adjacency must be changed. Thus $d_1$ is essentially the normalized $L_1$ distance between the two adjacency matrices. We say that $G$ is $\varepsilon$-far from being $F$-free if $d_1(G,H)\ge \varepsilon$ for every $F$-free graph $H$ on $[n]$. The graph removal lemma states that, if a graph is far from being $F$-free, then it should contain many copies of $F$. The precise statement is that, for every fixed graph $F$ and every $\varepsilon>0$, there is a constant $\delta_F(\varepsilon)>0$ such that every $n$-vertex graph that is $\varepsilon$-far from being $F$-free contains at least
$\delta_F(\varepsilon)n^{|V(F)|}$
copies of $F$. The triangle case is due to Ruzsa and Szemer\'edi~\cite{RuzsaSzemeredi}, and the general result was proved by Alon, Duke, Lefmann, R\"odl, and Yuster~\cite{AlonDukeLefmannRodlYuster} and by F\"uredi~\cite{Furedi}. The usual proof combines Szemer\'edi's regularity lemma~\cite{Szemeredi} with a counting lemma.

Removal lemmas are important not only in extremal combinatorics, but
also in number theory, discrete geometry, and theoretical computer
science. In particular, they give one-sided property testers. If every
graph that is $\varepsilon$-far with respect to the $L_1$ distance from being $F$-free contains many
copies of $F$, then a small random sample finds a copy with constant
probability. The quantitative bound on $\delta_F(\varepsilon)$
therefore controls the query complexity of such a tester.

A central question is to determine the best possible dependence of
$\delta_F(\varepsilon)$ on $\varepsilon$. This remains a major open
problem even for $F=K_3$. Alon~\cite{Alon} proved that a polynomial
bound is possible if and only if $F$ is bipartite, but this does not
determine the sharp quantitative bound. See
\cite{GishbolinerShapira} for a recent survey and for related results
for directed graphs, tournaments, and hypergraphs.

Permutations arise naturally in rankings and permutation patterns give a similar setting. A permutation $\pi$ contains a pattern $\sigma$ if some subsequence of $\pi$ has the same relative order as $\sigma$. Pattern avoidance is a basic topic in both extremal and enumerative combinatorics. This leads to the following property-testing question: can one decide whether a very large permutation is $\sigma$-free, or far from being $\sigma$-free, by looking at only a small random subpermutation?

For permutations, the answer depends on the metric. Cooper~\cite{Cooper} proved a permutation removal lemma using a regularity argument. However, the result was not stated with respect to a metric, so it did not directly give a property-testing result. Hoppen, Kohayakawa, Moreira, and Sampaio~\cite{HoppenEtAl} proved that hereditary permutation properties are one-sided testable under the rectangular distance. Their proof uses compactness and gives no effective bound on the query complexity.

Klimo\v{s}ov\'a and Kr\'al'~\cite{KlimosovaKral} proved that every hereditary permutation property is testable under Kendall's tau distance, the normalized number of pairs whose relative order differs. Their bound is of Ackermann type in $1/\varepsilon$.  Fox and Wei~\cite{FoxWeiPrep,FoxWeiSODA} obtained a polynomial bound for testing $\sigma$-freeness for every fixed pattern $\sigma$.

Fox and Wei~\cite{FoxWeiSODA,FoxWeiCPC} also studied the earth mover's distance
\[
\rho_{\mathrm{EMD}}(\pi_1,\pi_2)
 =\binom{n}{2}^{-1}
   \min_{\theta\in S_n}
   \sum_{i=1}^n
   \bigl(|i-\theta(i)|
        +|\pi_1(i)-\pi_2(\theta(i))|\bigr).
\]
This is the minimum normalized total $L_1$ transportation cost between
the point sets of the two permutations.

\begin{figure}[!htb]
\centering
\begin{minipage}[t]{0.48\textwidth}
 \vspace{0pt}
 \centering
 \resizebox{.7\linewidth}{!}{%
 \begin{tikzpicture}[scale=0.51]
  \draw[step=1,black!70,thin] (0,0) grid (8,8);
  \foreach \xa/\ya/\xb/\yb in {
   0/8/2/7,
   1/3/1/5,
   2/1/0/0,
   3/4/4/2,
   4/7/3/6,
   5/6/5/8,
   6/5/8/4,
   7/2/7/3,
   8/0/6/1}
  \draw[violet!85,dash pattern=on 4pt off 2pt,line width=0.8pt]
   (\xa,\ya)--(\xb,\yb);
  \foreach \x/\y in
   {0/0,1/5,2/7,3/6,4/2,5/8,6/1,7/3,8/4}
  \fill[red] (\x,\y) circle (2.5pt);
  \foreach \x/\y in
   {0/8,1/3,2/1,3/4,4/7,5/6,6/5,7/2,8/0}
  \fill[blue]
   (\x-0.09,\y-0.09) rectangle (\x+0.09,\y+0.09);
 \end{tikzpicture}%
 }
 \caption{An optimal matching (dashed line segments) with total $L_1$ cost $22$ from
 {\color{blue}$\pi_1=942587631$} to
 {\color{red}$\pi_2=168739245$}; hence
 $\rho_{\mathrm{EMD}}(\pi_1,\pi_2)=22/\binom{9}{2}=11/18$.}
 \label{fig:emd-example}
\end{minipage}\hfill%
\begin{minipage}[t]{0.48\textwidth}
 \vspace{0pt}
 \centering
 \resizebox{.7\linewidth}{!}{%
 \begin{tikzpicture}[scale=0.40]
  \draw[step=1,black!70,thin] (0,0) grid (14,14);

  \foreach \xa/\ya/\xb/\yb in {0/14/0/0,14/0/14/14}
  \draw[violet!85,dash pattern=on 4pt off 2pt,line width=1.6pt]
   (\xa,\ya)--(\xb,\yb);

  \foreach \x/\y in
   {0/14,1/13,2/12,3/11,4/10,5/9,6/8,7/7,8/6,9/5,10/4,11/3,12/2,13/1,14/0}
  \fill[blue]
   (\x-0.18,\y-0.18) rectangle (\x+0.18,\y+0.18);

  \foreach \x/\y in
   {0/0,1/13,2/12,3/11,4/10,5/9,6/8,7/7,8/6,9/5,10/4,11/3,12/2,13/1,14/14}
  \fill[red] (\x,\y) circle (3.5pt);
 \end{tikzpicture}%
 }
 \caption{An optimal matching ($\theta(i)=i$ for every $i\neq 1,n$) from the decreasing permutation
 {\color{blue}$\delta_{15}$} to its endpoint swap
 {\color{red}$\tau$} ($\tau(1)=1$ and $\tau (n)=n$) with total $L_1$ cost $28$; hence
 $\rho_{\mathrm{EMD}} (\delta_{15},\tau)=28/\binom{15}{2}=4/15$.}
 \label{fig:emd-endpoint-swap}
\end{minipage}
\end{figure}

Fox and
Wei~\cite{FoxWeiSODA,FoxWeiCPC} proved polynomial-query two-sided
testing results for hereditary permutation properties and also
obtained a canonical one-sided tester for a closely related metric
that is quantitatively equivalent to earth mover's distance. In
particular, for every fixed pattern $\sigma\in S_k$ and every
$\varepsilon>0$, their results imply a removal lemma with respect to earth mover's distance: There exists a
$\delta_{\sigma,\mathrm{EMD}}(\varepsilon)>0$ such that every
sufficiently large permutation that is $\varepsilon$-far from being
$\sigma$-free under earth mover's distance contains at least
$\delta_{\sigma,\mathrm{EMD}}(\varepsilon)n^k$
copies of $\sigma$.

The best possible dependence of
$\delta_{\sigma,\mathrm{EMD}}(\varepsilon)$ on $\varepsilon$ is not
known in general. Thus the sharp quantitative removal problem remains
open both for graph edit distance and for the permutation earth
mover's distance. Garbe, Hladk\'y, Kun, and
Pek\'arkov\'a~\cite{GarbeHladkyKunPekarkova} later gave a
nonquantitative proof of permutation removal using the structure of
pattern-avoiding permutons.

\subsection{Our result}

For the permutation metrics discussed above, farness from
$\sigma$-freeness forces a positive density of $\sigma$-copies. The
earth mover's distance is an $L_1$-type metric: it measures the total
transportation cost, so a small number of points may move a long
distance at small total cost. We study the corresponding bottleneck
metric, in which every point must move only a short distance.

For two permutations $\pi_1,\pi_2\in S_n$, define the {\it $L_\infty$-earth mover's distance} as 
\[
\rho_\infty(\pi_1,\pi_2)
 =\frac{1}{n-1}
   \min_{\theta\in S_n}
   \max_{1\le i\le n}
   \max\bigl\{|i-\theta(i)|,
               |\pi_1(i)-\pi_2(\theta(i))|\bigr\}.
\]

\begin{figure}[!htb]
\centering
\begin{minipage}[t]{0.48\textwidth}
 \vspace{0pt}
 \centering
 \resizebox{.7\linewidth}{!}{%
 \begin{tikzpicture}[scale=0.51]
  \draw[step=1,black!70,thin] (0,0) grid (8,8);
  \foreach \xa/\ya/\xb/\yb in {
   0/8/2/7,
   1/3/1/5,
   2/1/0/0,
   3/4/4/2,
   4/7/3/6,
   5/6/5/8,
   6/5/8/4,
   7/2/7/3,
   8/0/6/1}
  \draw[violet!85,dash pattern=on 4pt off 2pt,line width=0.8pt]
   (\xa,\ya)--(\xb,\yb);
  \foreach \x/\y in
   {0/0,1/5,2/7,3/6,4/2,5/8,6/1,7/3,8/4}
  \fill[red] (\x,\y) circle (2.5pt);
  \foreach \x/\y in
   {0/8,1/3,2/1,3/4,4/7,5/6,6/5,7/2,8/0}
  \fill[blue]
   (\x-0.09,\y-0.09) rectangle (\x+0.09,\y+0.09);
 \end{tikzpicture}%
 }
 \caption{An optimal matching (dashed line segments) with maximum $L_\infty$
 displacement $2$ from {\color{blue}$\pi_1=942587631$}
 to {\color{red}$\pi_2=168739245$}; hence
 $\rho_\infty(\pi_1,\pi_2)=2/(9-1)=1/4$.}
 \label{fig:rho-example}
\end{minipage}\hfill%
\begin{minipage}[t]{0.48\textwidth}
 \vspace{0pt}
 \centering
 \resizebox{.7\linewidth}{!}{%
 \begin{tikzpicture}[
  scale=0.40,
  midarrow/.style={
   violet!85,
   dash pattern=on 4pt off 2pt,
   line width=1pt,
   postaction={decorate},
   decoration={
    markings,
    mark=at position 0.5 with {
     \arrow{Stealth[length=3pt,width=3.5pt]}
    }
   }
  }
 ]
  \draw[step=1,black!70,thin] (0,0) grid (14,14);

  \foreach \xa/\ya/\xb/\yb in {
   0/14/1/13,
   1/13/2/12,
   2/12/3/11,
   3/11/4/10,
   4/10/5/9,
   5/9/6/8,
   6/8/0/0,
   7/7/14/14,
   8/6/7/7,
   9/5/8/6,
   10/4/9/5,
   11/3/10/4,
   12/2/11/3,
   13/1/12/2,
   14/0/13/1}
  \draw[midarrow] (\xa,\ya)--(\xb,\yb);

  \foreach \x in {0,...,14}
  \fill[blue]
   (\x-0.3,14-\x-0.3) rectangle
   (\x+0.3,14-\x+0.3);

  \foreach \x in {1,...,13}
  \fill[red] (\x,14-\x) circle (5pt);
  \fill[red] (0,0) circle (5pt);
  \fill[red] (14,14) circle (5pt);
 \end{tikzpicture}%
 }
 \caption{An optimal matching (dashed line segments) from the decreasing permutation
 {\color{blue}$\delta_{15}$} to its endpoint swap
 {\color{red}$\tau$} with maximum $L_\infty$ displacement $8$;
 hence $\rho_\infty(\delta_{15},\tau)=8/(15-1)=4/7$.}
 \label{fig:rho-endpoint-swap}
\end{minipage}
\end{figure}

Figure~\ref{fig:rho-example} and Figure \ref{fig:rho-endpoint-swap} provide a visual example of the distance between two permutations, illustrating an optimal geometric matching that minimizes the maximum individual displacement.
Thus earth mover's distance minimizes a sum of displacements, while
$\rho_\infty$ minimizes the largest displacement of any matched
point. 
Although $\rho_{\mathrm{EMD}}$ and $\rho_\infty$ have comparable
magnitudes for the illustrated pair $(\pi_1,\pi_2)$, the decreasing
permutation $\delta_n$ and its endpoint swap $\tau_n$ satisfy
$\rho_{\mathrm{EMD}}(\delta_n,\tau_n)=4/n$ while
$\rho_\infty(\delta_n,\tau_n)\to1/2$ as $n\to\infty$, showing that
$\rho_\infty$ can exceed $\rho_{\mathrm{EMD}}$ by an arbitrarily
large factor.
A small $\rho_\infty$ asks that every matched point move only a short distance.

We ask the following question: if a permutation is
$\varepsilon$-far from being $\sigma$-free under $\rho_\infty$, how
many copies of $\sigma$ must it contain?

\paragraph{Complete dichotomy.}
Let 
$\mathcal E
 =\{(12),(21),(132),(213),(231),(312)\}$ be the set of siz permutations.
We prove that the answer has only two possible forms. For
$\sigma\in\mathcal E$, $\varepsilon$-farness forces
$\Theta(\varepsilon n)$ copies. For every other pattern $\sigma$ of
length $k$, a single copy can keep a permutation $1/(8k)$-far from
being $\sigma$-free. Thus the usual positive-density removal
conclusion fails for every nontrivial pattern, but the precise sparse
behavior is completely classified.

To be more precise, for $\sigma\in S_k$ and $\pi\in S_n$, a \emph{$\sigma$-copy} in
$\pi$ is a set of indices $1\le i_1<\cdots<i_k\le n$ such that
\[
\pi(i_a)<\pi(i_b)
\quad\Longleftrightarrow\quad
\sigma(a)<\sigma(b)
\qquad\text{for all }a,b\in[k].
\]
We denote the number of $\sigma$-copies in $\pi$ by
$N_\sigma(\pi)$ and call $\pi$ \emph{$\sigma$-free} if
$N_\sigma(\pi)=0$.

\begin{theorem}\label{thm:linear}
Let $\sigma\in\mathcal E$. For every $n\ge2$, every
$\pi\in S_n$, and every $\varepsilon>0$, if $\pi$ is
$\varepsilon$-far from being $\sigma$-free with respect to
$\rho_\infty$, then
\[
N_\sigma(\pi)\ge \varepsilon(n-1).
\]
The bound $\varepsilon(n-1)$ is essentially tight. Indeed, for every $0<\varepsilon\le1/16$ and all sufficiently
large $n$, there is a permutation $\pi\in S_n$ that is
$\varepsilon$-far from being $\sigma$-free and satisfies
\[
N_\sigma(\pi)\le6\varepsilon n.
\]
\end{theorem}

\begin{theorem}
\label{thm:single}
Let $\sigma\in S_k$, where $k\ge2$ and
$\sigma\notin\mathcal E$. There exists $n_0(k)$ such that, for every
$n\ge n_0(k)$, there is a permutation $\pi\in S_n$ that is
$1/(8k)$-far from being $\sigma$-free with respect to
$\rho_\infty$, yet contains exactly one copy of $\sigma$.
\end{theorem}

Together, Theorems~\ref{thm:linear} and~\ref{thm:single} determine
the correct scale for every nontrivial pattern. They determine both
the power of $n$ and, in the six linear cases, the dependence on
$\varepsilon$ up to an absolute constant factor.

For example, for every sufficiently large $n$, there is a
permutation of length $n$ with exactly one increasing triple that is
$1/24$-far from having no increasing triple. By contrast, every
permutation that is $\varepsilon$-far from having no increasing pair
contains at least $\varepsilon(n-1)$ increasing pairs. Thus the
number of violations and the distance from the property can behave
in completely different ways.

Both sides of the dichotomy have the same consequence for canonical
one-sided property testing. In the examples from
Theorem~\ref{thm:linear}, the number of copies is $O(n)$, and in
Theorem~\ref{thm:single} it is one. In either case the density of
$\sigma$-copies tends to zero as $n\to\infty$. Therefore no
canonical one-sided tester that samples a bounded-size
subpermutation and rejects only after finding a $\sigma$-copy can
have query complexity depending only on $\varepsilon$ and $\sigma$.

\subsection{Proof ideas and challenges}

For graphs, the classical approach to removal lemmas combines regularity and counting lemmas. In the permutation setting, standard approaches include both regularity methods and compactness arguments: Cooper~\cite{Cooper} established permutation regularity and counting lemmas, while Garbe, Hladk\'y, Kun and Pek\'arkov\'a~\cite{GarbeHladkyKunPekarkova} gave a compactness-based proof of permutation removal. Fox, Lov\'asz and Zhao~\cite[Theorem~7.4] {FoxLovaszZhao} later improved the tower-type bound of Hoppen, Kohayakawa and Sampaio~\cite{HoppenKS} to a single-exponential bound on the number of parts in an equitable interval partition.

Our goal, however, is to determine the sharp dependence on both $\varepsilon$ and $n$. The compactness proof is nonquantitative, while the usual regularity-based removal proofs incur substantial quantitative losses in $\varepsilon$, in contrast to the linear dependence in Theorem~\ref{thm:linear}. More importantly, Cooper's counting theorem~\cite[Theorem~11]{Cooper} estimates the number of copies of a length-$k$ pattern up to an additive error of the form $\eta n^k$, where $\eta>0$ depends on the approximation parameters. With these parameters fixed independently of $n$, this estimate alone is too coarse to distinguish counts of order $\varepsilon n$ from a single copy or from none---precisely the scales relevant to Theorems~\ref{thm:linear} and~\ref{thm:single}.  Such small counts matter because $\rho_\infty$ takes the maximum
displacement over all points, while the earth mover's distance takes the
sum, so a few exceptional points cannot be absorbed into an error term.

For the lower bound in Theorem~\ref{thm:linear}, we destroy copies by a sequence of small rearrangements of $\pi$, each of which destroys at least as many copies as the largest horizontal or vertical distance by which it moves a point. Since $\pi$ is $\varepsilon$-far from being $\sigma$-free, the triangle inequality shows that these distances add up to at least $\varepsilon(n-1)$, and hence so does the number of destroyed copies. The key point is that a rearrangement must not create new copies. 
The case $\sigma=(12)$ is relatively simple, since the
$(12)$-avoiding permutations are precisely the decreasing
permutations. Swapping two consecutive elements can be such a rearrangement. 
For $\sigma=(132)$, however, destroying a chosen
copy by swapping entries can create new copies involving other
points, and the effect could propogate as we do more swappings. Moreover, changing only a few entries does not by itself
control $\rho_\infty$, which depends on the largest displacement.
We must therefore find a rearrangement that creates no new copies
and destroys enough existing copies to pay for its largest
displacement. We choose a $(132)$-copy whose first two entries
are adjacent and rearrange the prefix ending at its second entry,
placing entries above its third entry's value before those below
it while preserving order within each group. No resulting
$(132)$-copy can meet both groups, so every surviving copy pulls
back to an original copy. The same rearrangement destroys at
least $d$ copies while moving every point by at most $d$ in each
coordinate, for a suitable $d>0$; iteration and the triangle
inequality then give the lower bound in
Theorem~\ref{thm:linear}. 
The sharp upper bound constructions come from small modifications of the decreasing permutation.

For Theorem~\ref{thm:single}, we build $\pi$ explicitly so that it is both $\epsilon$-far from being $\sigma$-free but it only contains exactly one copy of $\sigma$. If the entries of a $\sigma$-copy are pairwise separated by at least $2\varepsilon(n-1)$ in both position and value, then this copy alone ensures that $\pi$ is $\varepsilon$-far from being $\sigma$-free under $\rho_\infty$. We therefore start with a single copy of $\sigma$, whose $|\sigma|$ points are well-separated among the $n$ points, and add the rest of the $n- |\sigma|$ points. Our goal is to add these $n- |\sigma|$  points into blocks of points, aiming to create the required separation between the original points in $\sigma$ while not creating any other copy of $\sigma$ using any of these $n- |\sigma|$ points. 


\begin{figure}[h]
\centering
\begin{tikzpicture}[scale=1.25,>=stealth]
  \draw[->,thin] (0,0) -- (4.3,0) node[right] {$x$};
  \draw[->,thin] (0,0) -- (0,4.3) node[above] {$y$};

  \draw[step=1,black!60,thin] (1,1) grid (3,3);
  \draw[black!60,thin] (1,1) rectangle (3,3);

  \foreach \x/\y/\lab/\pos in {
    1/1/{1}/below left,
    2/2/{2}/below left,
    3/3/{3}/above right}
  {
    \fill[blue] (\x,\y) circle (2.8pt);
    \node[\pos] at (\x,\y) {\small $\lab$};
  }

  \foreach \t in {0.12,0.24,...,0.96}
    \fill[blue] ({1-\t},{2+\t}) circle (1.2pt); 

  \foreach \t in {0.12,0.24,...,0.96}
    \fill[blue] ({2-\t},{3+\t}) circle (1.2pt); 

  \foreach \t in {0.12,0.24,...,0.96}
    \fill[blue] ({2+\t},{1-\t}) circle (1.2pt); 

  \foreach \t in {0.12,0.24,...,0.96}
    \fill[blue] ({3+\t},{2-\t}) circle (1.2pt); 

  \node[blue,left]        at (0.15,2.85) {$P_L$};
  \node[blue,above left]  at (1.15,3.85) {$P_U$};
  \node[blue,below right] at (2.85,0.15) {$P_B$};
  \node[blue,right]       at (3.95,1.15) {$P_R$};
\end{tikzpicture}
\caption{The permutation {\color{blue}$\pi$} of length $n$ is obtained by fixing a copy of $(123)$ and distributing the $n-3$ additional points as evenly as possible among $P_L,P_U,P_B,P_R$; then $\pi$ is $1/9$-far from being $(123)$-free and contains a unique $(123)$-copy.}
\label{fig:123-construction}
\end{figure}
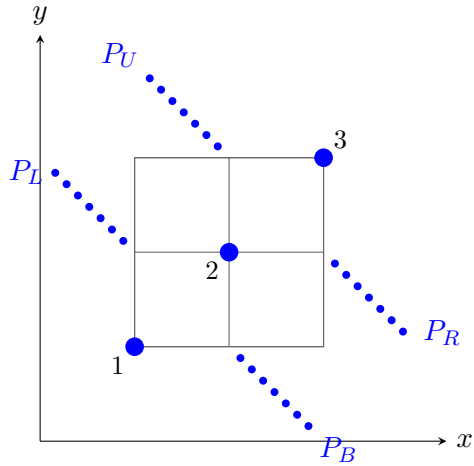

Figure~\ref{fig:123-construction}   illustrates this idea for $\sigma=(123)$. Distributing the $n-3$ additional points as evenly as possible along the lines $P_L,P_U,P_B,P_R$ ensures that, in the resulting permutation, the three consecutive entries of the original copy of $(123)$ are separated by at least $\lfloor(n-3)/4\rfloor$ in both position and value. 
For sufficiently large $n$, this is at least $2(n-1)/9$, so $\pi$ is $1/9$-far from being $(123)$-free. It is also straightforward to check that $\pi$ contains a unique $123$-copy with slight more effort. The main difficulty in proving Theorem~\ref{thm:single} for general $\sigma$ is to find constructions that achieve the required separation without creating any additional $\sigma$-copies.

Patterns $\sigma$ with little symmetry and large are also relatively simpler. Two entries $\sigma(i)$ and $\sigma(i+1)$ of a permutation $\sigma\in S_k$ are called \textit{twins} if $|\sigma(i+1)-\sigma(i)|=1$. A private initial conversation with Jacob Fox  suggests that, for sufficiently large twin-free patterns $\sigma\in S_k$, a short probabilistic construction may suffice: start with the canonical copy of $\sigma$ at the integer lattice points $(i,\sigma(i))$ of a $k\times k$ grid, and then insert blocks of additional points into suitably selected unit squares of the grid. However, unlike graphs, twins   are very likely in permutations. These locally similar entries make it harder to preserve the uniqueness of a $\sigma$-copy when adding blocks of points, even when the internal structure of each block is carefully designed. The main challenge is therefore to obtain constructions that work without a twin-free assumption, accommodate blocks containing twins of varying sizes and internal structures, and works for $\sigma$ of even small size. 


\begin{figure}[!htb]
\centering
\begin{minipage}[t]{0.48\textwidth}
\centering
\begin{tikzpicture}[scale=0.72,>=stealth]
  \fill[pink!35]   (0.5,0.5) rectangle (3.5,3.5);
  \fill[yellow!55] (3.5,3.5) rectangle (6.5,6.5);

  \draw[->,thin] (-3.5,0) -- (7.6,0) node[right] {$x$};
  \draw[->,thin] (0,-3.5) -- (0,7.6) node[above] {$y$};

  \draw[step=1,black!40,thin] (1,1) grid (6,6);
  \draw[black!70,thin] (1,1) rectangle (6,6);

  \foreach \xa/\ya/\xb/\yb in {
    1/2/0/3,
    0/3/-1/4,
    -1/4/-2/5,
    -2/5/-3/6,
    2/1/3/0,
    3/0/4/-1,
    4/-1/5/-2,
    5/-2/6/-3,
    2/6/1/7,
    6/2/7/1}
    \draw[blue,densely dotted,line width=0.9pt,
          shorten <=1pt,shorten >=1pt]
      (\xa,\ya)--(\xb,\yb);

  \node[blue,left]        at (-1.5,4.5) {$P_L$};
  \node[blue,above right] at (1.5,6.5)  {$P_U$};
  \node[blue,right]       at (4.5,-1.5) {$P_B$};
  \node[blue,above right] at (6.5,1.5)  {$P_R$};

  \foreach \x/\y in {1/1, 2/3,3/2,4/6,5/4,6/5}
    \fill[red] (\x,\y) circle (2.7pt);
\end{tikzpicture}
\caption{For $\sigma=(132645)$, adding any points on the blue segments
$P_L,P_U,P_B,P_R$ does not create a new copy of $(132645)$.}
\label{fig:132645-construction}
\end{minipage}\hfill
\begin{minipage}[t]{0.48\textwidth}
\centering
\begin{tikzpicture}[scale=0.72,>=stealth]
  \fill[pink!35]   (0.5,0.5) rectangle (1.5,1.5);
  \fill[yellow!55] (1.5,1.5) rectangle (6.5,6.5);

  \draw[->,thin] (-3.5,0) -- (7.6,0) node[right] {$x$};
  \draw[->,thin] (0,-3.5) -- (0,7.6) node[above] {$y$};

  \draw[step=1,black!40,thin] (1,1) grid (6,6);
  \draw[black!70,thin] (1,1) rectangle (6,6);

  \foreach \xa/\ya/\xb/\yb in {
    1/2/0/3,
    0/3/-1/4,
    -1/4/-2/5,
    -2/5/-3/6,
    2/1/3/0,
    3/0/4/-1,
    4/-1/5/-2,
    5/-2/6/-3,
    2/6/1/7,
    6/2/7/1}
    \draw[blue,densely dotted,line width=0.9pt,
          shorten <=1pt,shorten >=1pt]
      (\xa,\ya)--(\xb,\yb);

  \node[blue,left]        at (-1.5,4.5) {$P_L$};
  \node[blue,above right] at (1.5,6.5)  {$P_U$};
  \node[blue,right]       at (4.5,-1.5) {$P_B$};
  \node[blue,above right] at (6.5,1.5)  {$P_R$};

  \foreach \x/\y in {1/1,2/6,3/5,4/4,5/3,6/2}
    \fill[red] (\x,\y) circle (2.7pt);
\end{tikzpicture}
\caption{For $\sigma=(165432)$, the same construction fails: in fact,
any long decreasing point set in the region $x>1$, $y>1$ creates a new
copy of $(165432)$.}
\label{fig:165432-construction}
\end{minipage}
\end{figure}

The sizes and internal structures of such blocks also matter.
For example, both $(132645)$ and $(165432)$ are blow-ups of
$(12)$, each consists of two consecutive blocks
with every entry in the first block smaller than every entry
in the second, as illustrated in Figures~\ref{fig:132645-construction}
and~\ref{fig:165432-construction} where the colored squares
represent the two blocks. For $(132645)$, where both blocks
have three entries, the boundary-segment construction illustrated
in Figure~\ref{fig:132645-construction} is generalized from construction for $(123)$ in Figure \ref{fig:123-construction}, and is an instance of a
universal construction that works for all blow-ups of $(12)$
in which both blocks have size at least two and $\sigma(1)=1$. However, this
simple construction can fail when one  block in $\sigma$  has too few points, for example one
point, as in $(165432)$: adding even one point on $P_U$ or
$P_R$ in Figure~\ref{fig:165432-construction} creates another
$(165432)$-copy. The decreasing structure of the second block
presents a further difficulty: any five points forming a
decreasing sequence in the region $x>1$, $y>1$, together
with the point $(1,1)$, form another $(165432)$-copy.

We overcome these difficulties with a unified geometric approach
that yields constructions for all patterns covered by
Theorem~\ref{thm:single}, regardless of the sizes and internal structures
of their blocks.
We place most added segments outside the canonical square
$[1,k]\times[1,k]$, choosing their slopes and locations according
to how $\sigma$ begins and ends and where its smallest and
largest values occur. These choices ensure that each boundary
segment can contain only the corresponding first, last, lowest,
or highest entry of a competing copy. We combine these restrictions
with propagation arguments: a deviation from an original position
or value coordinate forces further deviations. The resulting
constraints force the competing copy to coincide with the original
one, unless $\sigma$ belongs to a special family that we treat
separately in Subsection \ref{subsec:two-blocks} and \ref{subsec:exceptional-five}.

When the boundary segments leave a required horizontal or vertical
strip unfilled, we add segments inside the grids. These internal
segments can interrupt the propagation, so we must control which
entries of a competing copy can lie on them and where the
propagation remains valid. In the most delicate case, we use two
internal segments and show that the existence of another copy
would force $\sigma=(31524)$, which is handled separately together
with its symmetric versions.


\paragraph{Organization of the paper.}
In Section~2, we introduce the necessary definitions and preliminaries, formalizing the geometric point-set representation of permutations and the $\rinf$ distance metric. In Section~3, we prove Theorem~\ref{thm:linear} by establishing linear lower bounds for the specific class of permutations that exhibit a removal-type phenomenon. Section~4 is dedicated to the proof of Theorem~\ref{thm:single}, where we detail the geometric constructions and structural lemmas required to embed a unique $\sigma$-copy while forcing a positive $\rinf$ distance from being $\sigma$-free.

\section{Definitions and preliminaries}

We write $[n]=\{1,2,\ldots,n\}$ and $S_n$ for the set of permutations of $[n]$.

In this section, we formalize the geometric representation of permutations, the specific infinity metric $\rinf$, and the separation criterion utilized throughout the paper.

To formally analyze permutations and their spatial displacements, it is convenient to represent them as point sets in $\mathbb{R}^2$. While any permutation $\pi$ of length $n$ can be naturally represented by lattice points in $[n]\times[n]$, proving our main geometric constructions requires relaxing this integer constraint.

\begin{definition}\label{def:associated-point-set}
A point set
\(\Pi=\{(x_1,y_1),(x_2,y_2),\ldots,(x_n,y_n)\}\subset\mathbb{R}^2\)
is associated with a permutation $\pi$ of length $n$ if:
\begin{itemize}[leftmargin=2em]
\item $x_1<x_2<\cdots<x_n$, and
\item the $y$-coordinates are pairwise distinct,
\end{itemize}
and $\pi(i)$ is the rank of $y_i$ among $\{y_1,\ldots,y_n\}$ when ordered in increasing order.

For such a point set, we denote by $\Pi(i)$ the unique point of $\Pi$ with $x$-coordinate $x_i$, and by $\Pi^{-1}(j)$ the unique point of $\Pi$ whose $y$-coordinate has rank $j$.
\end{definition}

As an immediate consequence of this ranking definition, for any point set $\Pi$ associated with a permutation $\pi$, we have $\Pi(i)=\Pi^{-1}(j)$ if and only if $\pi(i)=j$.

When a permutation $\pi\in S_n$ is represented by integer lattice points in $[n]\times[n]$, we refer to this as its canonical point set. Specifically, the canonical point set associated with $\pi$ is the integer lattice representation
\(\Pi_\pi=\{(i,\pi(i)):1\le i\le n\}\subset[n]\times[n]\).
In particular, $\Pi_\pi(i)=(i,\pi(i))$ and $\Pi_\pi^{-1}(j)=(\pi^{-1}(j),j)$.

With the point-set representation established, we can formalize the concept of spatial displacement. For any point $p\in\mathbb{R}^2$, let $p_x$ and $p_y$ denote its $x$- and $y$-coordinates, respectively. We define the standard $L_\infty$ distance between any two points $p,q\in\mathbb{R}^2$ by
\(\linf(p,q)=\max\{|p_x-q_x|,|p_y-q_y|\}\).

Using this geometry, we define the normalized infinity distance between permutations. This metric bounds the maximum $L_\infty$ distance any individual element must travel to map one permutation to another under an optimal matching. It is formally equivalent to a normalized $L_\infty$ bottleneck matching distance between the geometric point-set representations of the permutations.

\begin{definition}\label{def:rho-infty}
For any two permutations $\pi_1,\pi_2$ of length $n$, where $n\ge2$, the normalized infinity distance $\rinf$ is defined by
\(\rinf(\pi_1,\pi_2) =\frac{1}{n-1}\min\limits_{\theta\in S_n} \max\limits_{1\le i\le n} \max\bigl\{|i-\theta(i)|,|\pi_1(i)-\pi_2(\theta(i))|\bigr\}\).
\end{definition}

Let $\Pi_1$ and $\Pi_2$ be the canonical point sets associated with $\pi_1$ and $\pi_2$ on the grid $[n]\times[n]$. Notice that the inner maximum in the equation above is exactly the point-to-point distance $\linf(p,q)$, where $p=(i,\pi_1(i))\in\Pi_1$ and $q=(\theta(i),\pi_2(\theta(i)))\in\Pi_2$. Because any permutation $\theta\in S_n$ induces a bijection $\Phi:\Pi_1\to\Pi_2$ mapping $\Pi_1(i)$ to $\Pi_2(\theta(i))$, we can equivalently express the $\rinf$ distance as the normalized minimax distance between matched point sets:
\(\rinf(\pi_1,\pi_2) =\frac{1}{n-1} \min\limits_{\substack{\Phi:\Pi_1\to\Pi_2\\\Phi\text{ bijective}}} \max\limits_{p\in\Pi_1}\linf(p,\Phi(p))\).

\begin{remark}[Metric property]
The function $\rinf$ is a metric on $S_n$. Symmetry and definiteness are immediate from the matching formulation. For the triangle inequality, let $\Phi_{12}:\Pi_1\to\Pi_2$ and $\Phi_{23}:\Pi_2\to\Pi_3$ be optimal bijections. Their composition is a bijection from $\Pi_1$ to $\Pi_3$, and the triangle inequality for $L_\infty$ gives
\(\linf\bigl(p,(\Phi_{23}\circ\Phi_{12})(p)\bigr) \le \linf\bigl(p,\Phi_{12}(p)\bigr) + \linf\bigl(\Phi_{12}(p),\Phi_{23}(\Phi_{12}(p))\bigr)\).
Taking the maximum over $p$ and dividing by $n-1$ yields
\(\rinf(\pi_1,\pi_3) \le \rinf(\pi_1,\pi_2)+\rinf(\pi_2,\pi_3)\).
\end{remark}

With the $\rinf$ metric defined, we can now formalize what it means for a permutation to be far from a specific property. Let $\cP$ be a permutation property, and let $\cP_n$ denote the set of all permutations in $\cP$ of length $n$.

\begin{definition}\label{def:far}
For any $\eps>0$, a permutation $\pi$ of length $n$ is said to be $\eps$-far from $\cP$ if
\(\min\limits_{\gamma\in\cP_n}\rinf(\pi,\gamma)\ge\eps\).
\end{definition}

Intuitively, $\pi$ is $\eps$-far from a property $\cP$ if it is at least $\eps$-far (under $\rinf$) from every permutation of length $n$ satisfying $\cP$.

A natural way to force a permutation to be $\eps$-far from being $\sigma$-free is to embed a single copy of $\sigma$ whose elements are so widely separated that destroying the copy requires a large displacement.

\begin{fact}\label{fact:separation}
Suppose $\pi$ is a permutation of length $n$. If the subpermutation of $\pi$ induced by a set of indices $S=\{i_1,i_2,\ldots,i_k\}$ forms a $\sigma$-copy, and
\(\min\limits_{\substack{i,j\in S\\i\ne j}} \min\{|i-j|,|\pi(i)-\pi(j)|\} \ge 2\eps(n-1)\),
then $\pi$ is $\eps$-far from being $\sigma$-free.
\end{fact}

\begin{proof}
Suppose, for contradiction, that there is a $\sigma$-free permutation $\gamma\in S_n$ with $\rinf(\pi,\gamma)<\eps$. By Definition~\ref{def:rho-infty}, there is a bijection between the canonical point sets of $\pi$ and $\gamma$ under which every point moves by $L_\infty$ distance strictly less than $\eps(n-1)$. Write the image of $(i,\pi(i))$ as $(\theta(i),\gamma(\theta(i)))$.

Take distinct $i,j\in S$ with $i<j$. Since $j-i\ge2\eps(n-1)$, we have $\theta(j)-\theta(i)\ge (j-i)-|j-\theta(j)|-|i-\theta(i)|>2\eps(n-1)-2\eps(n-1)=0$, and hence $\theta(i)<\theta(j)$. Likewise, if $\pi(i)<\pi(j)$, then $\gamma(\theta(j))-\gamma(\theta(i))\ge \pi(j)-\pi(i)-|\pi(j)-\gamma(\theta(j))|-|\pi(i)-\gamma(\theta(i))|>0$.
Thus the matching preserves both the relative $x$-order and the relative $y$-order of all points indexed by $S$. Their images therefore form a $\sigma$-copy in $\gamma$, contradicting that $\gamma$ is $\sigma$-free. Hence every $\sigma$-free $\gamma$ has $\rinf(\pi,\gamma)\ge\eps$.
\end{proof}

This fact provides a sufficient geometric condition for being $\eps$-far. If the elements forming the $\sigma$-copy are separated so that both coordinate gaps between every pair are at least $2\eps(n-1)$, eliminating this copy requires altering the relative spatial ordering of at least two of these points. Because they are separated in each direction by an initial gap of at least $2\eps(n-1)$, any matching with a different permutation that inverts their relative ordering in either direction must bridge this distance. By the triangle inequality, at least one of the points must travel an $L_\infty$ distance of at least $\eps(n-1)$, which corresponds to a normalized $\rinf$ displacement of $\eps$.

\section[O(n) cases]{$O(n)$ cases}\label{sec:linear-cases}

For permutations $\sigma\in S_k$ and $\pi\in S_n$, let $N_\sigma(\pi)$ denote the number of $\sigma$-copies contained in $\pi$. Let
\(\Pi=\{(i,\pi(i)):i\in[n]\}\subseteq[n]\times[n]\)
be the canonical point set associated with $\pi$. For $i,j\in[n]$, Definition~2.1 gives
\(\Pi(i)=(i,\pi(i))\) and \(\Pi^{-1}(j)=(\pi^{-1}(j),j)\).

The primary ingredient in the proof of Theorem~1.1 is the following lemma.

\begin{lemma}\label{lem:reduction}
If
\(\sigma\in\{(12),(21),(132),(231),(312),(213)\}\)
and $\pi\in S_n$ satisfies $N_\sigma(\pi)>0$, then there exists a distinct permutation $\pi'\in S_n\setminus\{\pi\}$ such that
\(N_\sigma(\pi)-N_\sigma(\pi')\ge (n-1)\rho_\infty(\pi,\pi')\).
\end{lemma}

To prove Lemma~\ref{lem:reduction}, we treat two representative patterns, $(12)$ and $(132)$. In each case, we construct $\pi'$ by defining a bijection $\Phi$ from the canonical point set $\Pi$ of $\pi$ to the canonical point set $\Pi'$ of $\pi'$. We prove that every $\sigma$-copy in $\Pi'$ has a $\sigma$-copy as its preimage under $\Phi$, so the construction creates no new copies. We then exhibit sufficiently many copies in $\Pi$ whose images are not copies in $\Pi'$ and compare their number with the maximum displacement under $\Phi$.

Iterating Lemma~\ref{lem:reduction} strictly decreases the number of $\sigma$-copies $N_\sigma(\pi)$ until a $\sigma$-free permutation is reached. The triangle inequality for $\rho_\infty$ then converts the sum of the stepwise $\rho_\infty$ into the lower bound in Theorem~1.1. At the end of the section, we give explicit constructions showing that the resulting linear bound is asymptotically tight up to a constant factor.

We will repeatedly use the dihedral symmetries of a permutation plot. For $\pi\in S_n$, reversal, complementation, and inversion are defined by
\(\pi^{\mathrm r}(i)=\pi(n+1-i)\), \(\pi^{\mathrm c}(i)=n+1-\pi(i)\), \(\pi^{-1}\).
On canonical point sets, these operations are respectively reflection across a vertical line, reflection across a horizontal line, and reflection across the main diagonal. Each reflection is an $L_\infty$-isometry. Consequently, for every symmetry $T$ generated by these operations,
\(\rho_\infty(T\pi,T\pi')=\rho_\infty(\pi,\pi')\) and \(N_{T\sigma}(T\pi)=N_\sigma(\pi)\),
where $T\sigma$ denotes the pattern obtained by applying the same symmetry to the plot of $\sigma$. Thus any bound or construction proved for $\sigma$ transfers to every pattern in its symmetry orbit.

\subsection[Patterns 12 and 21]{$\sigma\in\{(12),(21)\}$}

Because $(21)$ is the reverse (equivalently, the complement) of $(12)$, the symmetry identities above show that it suffices to prove Lemma~\ref{lem:reduction} for $\sigma=(12)$.

\begin{claim}\label{clm:adjacent-12}
For every $\pi\in S_n$ with $N_{(12)}(\pi)>0$, there exists a $(12)$-copy whose corresponding point set $\Sigma$ satisfies
\(\Sigma(2)_x-\Sigma(1)_x=1\).
\end{claim}

\begin{proof}
Among all point sets $\Sigma\subseteq\Pi$ associated with $(12)$-copies, choose $\Sigma^*$ minimizing the horizontal gap $\Sigma(2)_x-\Sigma(1)_x$.

Write $\Sigma^*=\{\Pi(i),\Pi(j)\}$, where $1\le i<j\le n$ and $\pi(i)<\pi(j)$. Suppose for contradiction that $j-i\ge2$. If $\pi(i+1)>\pi(i)$, then $\{\Pi(i),\Pi(i+1)\}$ is a $(12)$-copy with smaller horizontal gap. If $\pi(i+1)<\pi(i)$, then
\(\pi(i+1)<\pi(i)<\pi(j)\),
so $\{\Pi(i+1),\Pi(j)\}$ is a $(12)$-copy with horizontal gap $j-i-1$. Either case contradicts the minimality of $\Sigma^*$. Hence $j=i+1$.
\end{proof}

\begin{proof}[Proof of Lemma~\ref{lem:reduction} for $\sigma=(12)$]
Let
\(\Sigma^*=\{\Pi(i^*),\Pi(i^*+1)\}\)
be the copy supplied by Claim~\ref{clm:adjacent-12}. Thus $\pi(i^*)<\pi(i^*+1)$.

\begin{construction}\label{con:swap}

Define $\Phi:\Pi\to[n]\times[n]$ by setting $\Phi(\Pi(i^*))=(i^*+1,\pi(i^*))$, $\Phi(\Pi(i^*+1))=(i^*,\pi(i^*+1))$, and $\Phi(p)=p$ for every other $p\in\Pi$. Let $\Pi'=\Phi(\Pi)$. Then $\Pi'$ is the canonical point set of the permutation $\pi'\in S_n$ given by $\pi'(i^*)=\pi(i^*+1)$, $\pi'(i^*+1)=\pi(i^*)$, and $\pi'(i)=\pi(i)$ for $i\notin\{i^*,i^*+1\}$.
\end{construction}

Construction~\ref{con:swap} preserves every $y$-coordinate, swaps the two adjacent $x$-coordinates $i^*$ and $i^*+1$, and fixes every other $x$-coordinate. Hence
\inlinetag{\Phi(\Pi(a))_x<\Phi(\Pi(b))_x\text{ and }\{a,b\}\ne\{i^*,i^*+1\}\quad\Rightarrow\quad a<b}{1}{eq:swap-monotonicity}.
Moreover, the image $\Phi(\Sigma^*)$ is a $(21)$-copy, not a $(12)$-copy.

\begin{claim}\label{clm:preimage-12}
If $\Sigma'\subseteq\Pi'$ is associated with $(12)$, then $\Phi^{-1}(\Sigma')\subseteq\Pi$ is also associated with $(12)$.
\end{claim}

\begin{proof}
Write
\(\Sigma'=\{\Phi(\Pi(a)),\Phi(\Pi(b))\}\)
with
\(\Phi(\Pi(a))_x<\Phi(\Pi(b))_x\) and \(\Phi(\Pi(a))_y<\Phi(\Pi(b))_y\).
The pair cannot be $\{\Phi(\Pi(i^*)),\Phi(\Pi(i^*+1))\}$, because that pair is a $(21)$-copy. Thus \eqref{eq:swap-monotonicity} gives $a<b$. Since $\Phi$ preserves $y$-coordinates, $\pi(a)<\pi(b)$. Therefore $\{\Pi(a),\Pi(b)\}$ is a $(12)$-copy in $\Pi$.
\end{proof}

Claim~\ref{clm:preimage-12} gives an injection from the $(12)$-copies of $\pi'$ to those of $\pi$. The copy $\Sigma^*$ is not in the image of this injection, because $\Phi(\Sigma^*)$ is not a $(12)$-copy. Hence
\(N_{(12)}(\pi)-N_{(12)}(\pi')\ge1\).

For $i=i^*$ or $i=i^*+1$, the map $\Phi$ changes the $x$-coordinate by exactly one and leaves the $y$-coordinate unchanged, so
\(L_\infty(\Pi(i),\Phi(\Pi(i)))=1\).

For all other $i$, $\Phi(\Pi(i))=\Pi(i)$, so the point displacement is zero. Using $\Phi$ as a matching between the two canonical point sets therefore gives
\(\rho_\infty(\pi,\pi') \le \frac{1}{n-1}\max\limits_{1\le i\le n}L_\infty(\Pi(i),\Phi(\Pi(i))) =\frac{1}{n-1}\).

Combining the two estimates gives
\(N_{(12)}(\pi)-N_{(12)}(\pi') \ge1 \ge(n-1)\rho_\infty(\pi,\pi')\),
which proves the lemma for $(12)$ and, by symmetry, for $(21)$.
\end{proof}

\subsection[Patterns 132, 231, 312, and 213]{$\sigma\in\{(132),(231),(312),(213)\}$}

The reverse, complement, and reverse-complement of $(132)$ are $(231)$, $(312)$, and $(213)$, respectively. Therefore the symmetry identities above reduce all four cases to $\sigma=(132)$.

For a $(132)$-copy, write its associated point set as
\(\Sigma=\{(x_1,y_1),(x_2,y_2),(x_3,y_3)\}\), \(x_1<x_2<x_3\), \(y_1<y_3<y_2\).
Thus
\(\Sigma(1)=\Sigma^{-1}(1)=(x_1,y_1)\), \(\Sigma(2)=\Sigma^{-1}(3)=(x_2,y_2)\), \(\Sigma(3)=\Sigma^{-1}(2)=(x_3,y_3)\).

\begin{claim}\label{clm:adjacent-132}
For every $\pi\in S_n$ with $N_{(132)}(\pi)>0$, there exists a $(132)$-copy whose corresponding point set $\Sigma$ satisfies
\(\Sigma(2)_x-\Sigma(1)_x=1\).
\end{claim}

\begin{proof}
Among all point sets $\Sigma\subseteq\Pi$ associated with $(132)$-copies, choose $\Sigma^*$ minimizing the horizontal gap $\Sigma(2)_x-\Sigma(1)_x$.

Write $\Sigma^*=\{\Pi(i),\Pi(j),\Pi(k)\}$ with $1\le i<j<k\le n$ and
\(\pi(i)<\pi(k)<\pi(j)\).
Suppose for contradiction that $j-i\ge2$. If $\pi(i+1)<\pi(k)$, then $\{\Pi(i+1),\Pi(j),\Pi(k)\}$ is a $(132)$-copy. If $\pi(i+1)>\pi(k)$, then $\{\Pi(i),\Pi(i+1),\Pi(k)\}$ is a $(132)$-copy. In either case, the gap between the first two indices is at most
\(\max\{j-i-1,1\}<j-i\),
contradicting the minimality of $\Sigma^*$. Hence $j=i+1$.
\end{proof}

\begin{proof}[Proof of Lemma~\ref{lem:reduction} for $\sigma=(132)$]
Let
\(\Sigma^*=\{\Pi(x_1),\Pi(x_2),\Pi(x_3)\}\), \(1\le x_1<x_2<x_3\le n\),
be the copy supplied by Claim~\ref{clm:adjacent-132}.

In particular,
\(x_2=x_1+1\) and \(\pi(x_1)<\pi(x_3)<\pi(x_2)\).

\begin{construction}\label{con:stable-partition}
Partition the points of $\Pi$ into $U=\{(x,y):x\le x_2,\ y>\pi(x_3)\}$, $L=\{(x,y):x\le x_2,\ y<\pi(x_3)\}$, and $R=\{(x,y):x>x_2\}$. The omitted strip $\{(x,y):x\le x_2,\ y=\pi(x_3)\}$ contains no point of $\Pi$.

Let $U_x=\{i_1<i_2<\cdots<i_u=x_2\}$ and $L_x=\{j_1<j_2<\cdots<j_\ell\}$ be the sets of $x$-coordinates of $\Pi\cap U$ and $\Pi\cap L$, respectively. Define $\Phi:\Pi\to[n]\times[n]$ by $\Phi(\Pi(i_s))=(s,\pi(i_s))$ for $1\le s\le u$, by $\Phi(\Pi(j_t))=(u+t,\pi(j_t))$ for $1\le t\le\ell$, and by $\Phi(p)=p$ for $p\in R$. Let $\Pi'=\Phi(\Pi)$. Then $\Pi'$ is the canonical point set of the permutation $\pi'\in S_n$ determined by $\pi'(s)=\pi(i_s)$ for $1\le s\le u$, $\pi'(u+t)=\pi(j_t)$ for $1\le t\le\ell$, and $\pi'(r)=\pi(r)$ for $r>x_2$.
\end{construction}

Since the first $x_2$ points lie in exactly one of $U$ and $L$, we have $u+\ell=x_2$. Moreover, $x_1\in L_x$ and $x_2\in U_x$, so $u,\ell\ge1$.

Construction~\ref{con:stable-partition} preserves every $y$-coordinate. Among the first $x_2$ positions, it stably moves all points of $U$ to the left of all points of $L$; it fixes $R$ pointwise. Since $x_1=x_2-1$ lies in $L_x$ whereas $x_2$ lies in $U_x$, their relative $x$-order is reversed, so $\pi'\ne\pi$.

\begin{fact}\label{fact:stable-properties}
The bijection $\Phi$ satisfies the following properties.
\begin{enumerate}[label=\arabic*.,leftmargin=2em]
\item $\Phi(\Pi\cap U)\subseteq U$, $\Phi(\Pi\cap L)\subseteq L$, and $\Phi(\Pi\cap R)\subseteq R$.
\item If either $p,q\in U\cup R$ or $p,q\in L\cup R$, then
\(p_x<q_x\quad\Leftrightarrow\quad \Phi(p)_x<\Phi(q)_x\).
\item If $i<j_1$, then $\Phi(\Pi(i))=\Pi(i)$.
\end{enumerate}
\end{fact}

\begin{proof}
Property 1 follows directly from the coordinate ranges in Construction~\ref{con:stable-partition}. For Property 2, the equivalence is immediate when both points lie in $R$, since $\Phi$ fixes $R$ pointwise.

For the remaining pairs, the construction preserves the order inside each of $U$ and $L$. It also places every image point of $U$ or $L$ at an $x$-coordinate at most $x_2$, while every point of $R$ has $x$-coordinate greater than $x_2$. Hence the displayed equivalence holds both within a block and between that block and $R$.

For Property 3, the minimality of $j_1$ implies that every index $1,\ldots,j_1-1$ belongs to $U_x$. Thus $i_s=s$ for $s<j_1$, and the definition of $\Phi$ fixes $\Pi(s)$ for every such $s$.
\end{proof}

\begin{fact}\label{fact:no-mixed-copy}
No $(132)$-copy in $\Pi'$ intersects both $U$ and $L$.
\end{fact}

\begin{proof}
Every point of $\Pi'\cap U$ has $x$-coordinate at most $u$ and $y$-coordinate greater than $\pi(x_3)$, whereas every point of $\Pi'\cap L$ has $x$-coordinate at least $u+1$ and $y$-coordinate less than $\pi(x_3)$.

Suppose $\Sigma'\subseteq\Pi'$ is a $(132)$-copy. If its leftmost point $\Sigma'(1)$ lies in $U$, then
\(\Sigma'(1)_y<\Sigma'(3)_y<\Sigma'(2)_y\)
forces all three points to have $y$-coordinate greater than $\pi(x_3)$, so none lies in $L$. If $\Sigma'(1)$ lies in $L\cup R$, then $\Sigma'(1)_x>u$, and the two remaining points lie still farther to the right; hence none lies in $U$. Thus $\Sigma'$ cannot meet both blocks.
\end{proof}

\begin{claim}\label{clm:preimage-132}
If $\Sigma'\subseteq\Pi'$ is associated with $(132)$, then $\Phi^{-1}(\Sigma')\subseteq\Pi$ is also associated with $(132)$.
\end{claim}

\begin{proof}
By Fact~\ref{fact:no-mixed-copy}, $\Sigma'$ is contained in either $U\cup R$ or $L\cup R$.

On either of these two unions, Fact~\ref{fact:stable-properties} shows that $\Phi$ preserves the complete $x$-order, and Construction~\ref{con:stable-partition} preserves every $y$-coordinate. Therefore $\Phi^{-1}(\Sigma')$ has the same relative $x$- and $y$-orders as $\Sigma'$, and hence is a $(132)$-copy.
\end{proof}

We next identify $x_2-j_1$ distinct $(132)$-copies in $\Pi$ whose images are not copies. For each integer $x$ with $j_1+1\le x\le x_2$, define $C_x=\{\Pi(j_1),\Pi(x),\Pi(x_3)\}$ when $x\in U_x$, and define $C_x=\{\Pi(x),\Pi(x_2),\Pi(x_3)\}$ when $x\in L_x$.
If $x\in U_x$, then
\(j_1<x\le x_2<x_3\) and \(\pi(j_1)<\pi(x_3)<\pi(x)\).
If $x\in L_x$, then
\(x<x_2<x_3\) and \(\pi(x)<\pi(x_3)<\pi(x_2)\).
Thus every $C_x$ is a $(132)$-copy. The sets $C_x$ are distinct, and each intersects both $U$ and $L$. By Fact~\ref{fact:no-mixed-copy}, $\Phi(C_x)$ is not a $(132)$-copy.

Claim~\ref{clm:preimage-132} gives an injection from the $(132)$-copies of $\pi'$ to those of $\pi$, while the $x_2-j_1$ copies above are not in its image. Therefore
\inlinetag{N_{(132)}(\pi)-N_{(132)}(\pi')\ge x_2-j_1}{2}{eq:destroyed-132}.
We now bound the distance between $\pi$ and $\pi'$ under $\rho_\infty$. By the matching definition,
\(\rho_\infty(\pi,\pi') \le\frac{1}{n-1}\max\limits_{1\le i\le n}L_\infty(\Pi(i),\Phi(\Pi(i)))\).

Since $\Phi$ preserves $y$-coordinates,
\(L_\infty(\Pi(i),\Phi(\Pi(i))) =\bigl|\Pi(i)_x-\Phi(\Pi(i))_x\bigr|\).
If $\Pi(i)\in R$, this displacement is zero. If $\Pi(i)\in U$ and $i<j_1$, it is also zero by Fact~\ref{fact:stable-properties}. If $\Pi(i)\in U$ and $i>j_1$, then all indices $1,\ldots,j_1-1$ belong to $U_x$, so the rank of $i$ in $U_x$ lies in $[j_1,x_2]$; both its old and new $x$-coordinates therefore lie in $[j_1,x_2]$.

Finally, if $\Pi(i)\in L$, then $i\ge j_1$, while its new coordinate $u+t$ lies in $[j_1,x_2]$ because $u\ge j_1-1$.

Hence in every case
\(L_\infty(\Pi(i),\Phi(\Pi(i)))\le x_2-j_1\).

Using $\Phi$ as a matching gives
\(\rho_\infty(\pi,\pi') \le\frac{x_2-j_1}{n-1}\).
Together with \eqref{eq:destroyed-132}, this yields
\(N_{(132)}(\pi)-N_{(132)}(\pi') \ge(n-1)\rho_\infty(\pi,\pi')\),
which proves Lemma~\ref{lem:reduction} for $(132)$ and hence, by symmetry, for $(231)$, $(312)$, and $(213)$.
\end{proof}

\subsection{Proof of Theorem~1.1}

\begin{proof}[Proof of the lower bound in Theorem~1.1]
Starting with $\pi_0=\pi$, apply Lemma~\ref{lem:reduction} whenever the current permutation contains a $\sigma$-copy. The lemma produces $\pi_{i+1}\ne\pi_i$; since $\rho_\infty$ is a metric, its right-hand side is positive, and hence the integer $N_\sigma(\pi_i)$ strictly decreases. Therefore the process terminates at a $\sigma$-free permutation $\pi_m$. For every $0\le i<m$,
\(N_\sigma(\pi_i)-N_\sigma(\pi_{i+1}) \ge(n-1)\rho_\infty(\pi_i,\pi_{i+1})\).
Summing and using the triangle inequality gives $N_\sigma(\pi)=N_\sigma(\pi_0)-N_\sigma(\pi_m)\ge(n-1)\sum\limits_{i=0}^{m-1}\rho_\infty(\pi_i,\pi_{i+1})\ge(n-1)\rho_\infty(\pi_0,\pi_m)$.
If $\pi$ is $\varepsilon$-far from being $\sigma$-free, then $\rho_\infty(\pi_0,\pi_m)\ge\varepsilon$, because $\pi_m$ is $\sigma$-free. Therefore
\(N_\sigma(\pi)\ge\varepsilon(n-1)\).
\end{proof}

\subsection[Upper-bound constructions with O(epsilon n) copies]{Upper-bound constructions with $O(\varepsilon n)$ copies}

Fix
\(0<\varepsilon\le\frac1{16}\) and \(n\ge\frac1\varepsilon\),
and set
\(d=\left\lceil2\varepsilon(n-1)\right\rceil\).
Since $n\ge \frac{1}{\varepsilon}$ and $\varepsilon\le1/16$, we have $2\varepsilon(n-1)\ge2(1-\varepsilon)>1$, and hence $d\ge2$. Moreover,
\(d\ge2\varepsilon(n-1)\), \(d\le2\varepsilon(n-1)+1\le3\varepsilon n\),
and hence
\(4d-1\le12\varepsilon n-1\le\frac34n-1<n\).
These inequalities ensure that both constructions below are well-defined and that their exact copy counts are $O(\varepsilon n)$. The constructions for the symmetric patterns are obtained by applying the corresponding plot symmetries.

\subsubsection[Patterns 12 and 21]{$\sigma\in\{(12),(21)\}$}

By symmetry, it suffices to construct the example for $\sigma=(12)$.

\begin{construction}\label{con:upper-12}
Define $\pi\in S_n$ by $\pi(1)=n-d$, $\pi(d+1)=n$, and $\pi(i)=n-i+1$ for every $i\notin\{1,d+1\}$.
Equivalently, $\pi$ is obtained from the decreasing permutation by swapping the values at positions $1$ and $d+1$.
\end{construction}

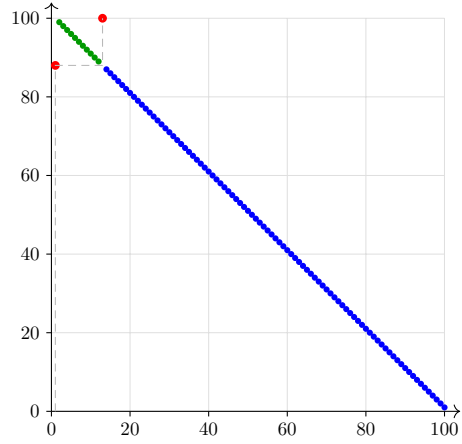
\begin{figure}[H]
\centering
\begin{tikzpicture}[x=0.052cm,y=0.052cm]
 \draw[step=20,gray!25,very thin] (0,0) grid (100,100);
 \draw[->] (0,0) -- (104,0);
 \draw[->] (0,0) -- (0,104);
 \foreach \t in {0,20,40,60,80,100}{
 \draw (\t,0) -- (\t,-1.3) node[below,scale=.65] {\t};
 \draw (0,\t) -- (-1.3,\t) node[left,scale=.65] {\t};
 }
 \foreach \i in {14,...,100}{
 \pgfmathtruncatemacro{\y}{101-\i}
 \fill[blue] (\i,\y) circle[radius=1.15pt];
 }
 \foreach \i in {2,...,12}{
 \pgfmathtruncatemacro{\y}{101-\i}
 \fill[green!60!black] (\i,\y) circle[radius=1.15pt];
 }
 \fill[red] (1,88) circle[radius=1.65pt];
 \fill[red] (13,100) circle[radius=1.65pt];
 \draw[densely dashed,gray!55] (1,0)--(1,88)--(13,88)--(13,100);
\end{tikzpicture}
\caption{$\pi$ from Construction~\ref{con:upper-12} for $n=100$ and $d=12$. The red points form the separated copy indexed by $S=\{1,d+1\}$; the green points participate in the remaining $(12)$-copies.}
\label{fig:upper-12}
\end{figure}

The indices $S=\{1,d+1\}$ induce a $(12)$-copy. Their horizontal and vertical separations are both $d$, because
\((d+1)-1=d\) and \(\pi(d+1)-\pi(1)=n-(n-d)=d\).
Since $d\ge2\varepsilon(n-1)$, Fact~2.4 implies that $\pi$ is $\varepsilon$-far from being $(12)$-free.

All entries other than those at positions $1$ and $d+1$ occur in decreasing order. The point at position $1$ forms a $(12)$-copy with each of the $d$ positions $2,\ldots,d+1$, and the point at position $d+1$ forms an additional $(12)$-copy with each of the $d-1$ positions $2,\ldots,d$. There are no other $(12)$-copies. Hence
\(N_{(12)}(\pi)=d+(d-1)=2d-1\le6\varepsilon n\).
Applying a reverse or complement gives the corresponding construction for $(21)$.

\subsubsection[Patterns 132, 231, 312, and 213]{$\sigma\in\{(132),(231),(312),(213)\}$}

By symmetry, it suffices to construct the example for $\sigma=(132)$.

\begin{construction}\label{con:upper-132}
Define $\pi\in S_n$ by
\[
\pi(i)=
\begin{cases}
n,&i=3d-1,\\
n-d,&i=4d-1,\\
n-i,&1\le i\le d-1\text{ or }3d\le i\le4d-2,\\
n-i-1,&d\le i\le3d-2,\\
n-i+1,&i\ge4d.
\end{cases}
\]
\end{construction}

The clauses in Construction~\ref{con:upper-132} use the six pairwise disjoint value sets $\{1,\ldots,n-4d+1\}$, $\{n-4d+2,\ldots,n-3d\}$, $\{n-3d+1,\ldots,n-d-1\}$, $\{n-d\}$, $\{n-d+1,\ldots,n-1\}$, and $\{n\}$. Their union is $[n]$, so the displayed formula indeed defines a permutation.

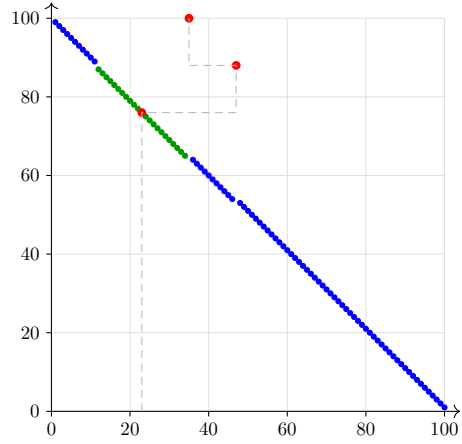
\begin{figure}[H]
\centering
\begin{tikzpicture}[x=0.052cm,y=0.052cm]
 \draw[step=20,gray!25,very thin] (0,0) grid (100,100);
 \draw[->] (0,0) -- (104,0);
 \draw[->] (0,0) -- (0,104);
 \foreach \t in {0,20,40,60,80,100}{
 \draw (\t,0) -- (\t,-1.3) node[below,scale=.65] {\t};
 \draw (0,\t) -- (-1.3,\t) node[left,scale=.65] {\t};
 }
 \foreach \i in {1,...,11}{
 \pgfmathtruncatemacro{\y}{100-\i}
 \fill[blue] (\i,\y) circle[radius=1.15pt];
 }
 \foreach \i in {12,...,22}{
 \pgfmathtruncatemacro{\y}{99-\i}
 \fill[green!60!black] (\i,\y) circle[radius=1.15pt];
 }
 \foreach \i in {24,...,34}{
 \pgfmathtruncatemacro{\y}{99-\i}
 \fill[green!60!black] (\i,\y) circle[radius=1.15pt];
 }
 \foreach \i in {36,...,46}{
 \pgfmathtruncatemacro{\y}{100-\i}
 \fill[blue] (\i,\y) circle[radius=1.15pt];
 }
 \foreach \i in {48,...,100}{
 \pgfmathtruncatemacro{\y}{101-\i}
 \fill[blue] (\i,\y) circle[radius=1.15pt];
 }
 \fill[red] (23,76) circle[radius=1.65pt];
 \fill[red] (35,100) circle[radius=1.65pt];
 \fill[red] (47,88) circle[radius=1.65pt];
 \draw[densely dashed,gray!55] (23,0)--(23,76)--(47,76)--(47,88)--(35,88)--(35,100);
\end{tikzpicture}
\caption{$\pi$ from Construction~\ref{con:upper-132} for $n=100$ and $d=12$. The red points form the separated copy indexed by $S=\{2d-1,3d-1,4d-1\}$; the green points mark the other possible first elements of the $(132)$-copies.}
\label{fig:upper-132}
\end{figure}

The indices
\(S=\{2d-1,3d-1,4d-1\}\)
induce a $(132)$-copy, because
\(\pi(2d-1)=n-2d\), \(\pi(3d-1)=n\), \(\pi(4d-1)=n-d\).
Every pair of these three points is separated by at least $d$ in both coordinates. Fact~2.4 therefore implies that $\pi$ is $\varepsilon$-far from being $(132)$-free.

For every $j\notin\{3d-1,4d-1\}$ and every $i<j$, the piecewise definition gives $\pi(i)>\pi(j)$. Now let $a<b<c$ form a $(132)$-copy, so
\(\pi(a)<\pi(c)<\pi(b)\).
The index $c$ cannot be ordinary, and it cannot equal $3d-1$ because $\pi(3d-1)=n$. Hence $c=4d-1$. To have $\pi(b)>n-d$, either $b\in\{1,\ldots,d-1\}$ or $b=3d-1$. The first alternative is impossible because every earlier value is then larger than $\pi(b)$, so necessarily $b=3d-1$. Finally, since $a<b=3d-1$,
\(\pi(a)<n-d \quad\Leftrightarrow\quad d\le a\le3d-2\).
There are exactly $2d-1$ such choices. Consequently,
\(N_{(132)}(\pi)=2d-1\le6\varepsilon n\).
The remaining three patterns follow by reverse, complement, and reverse-complement symmetry.

\section[O(1) cases]{$O(1)$ cases}\label{sec:constant-cases}

The proof of Theorem~1.2 rests on the following geometric construction theorem. Throughout this section, $k$ denotes the length of the forbidden pattern and $N$ denotes the length of the host permutation.

\begin{theorem}\label{thm:geometric-construction}
Let $\sigma\in S_k$, where either $k\ge4$ or $\sigma\in\{(123),(321)\}$. There exists $N_0(k)$ such that, for every $N\ge N_0(k)$, there is an $N$-point set $\Pi\subset\RR^2$ associated with a permutation $\pi\in S_N$ with the following properties. If
\(\Sigma=\{(i,\sigma(i)):i\in[k]\}\)
is the canonical point set of $\sigma$, then $\Sigma\subseteq\Pi$ and $\Sigma$ is the unique $\sigma$-copy in $\Pi$. Moreover, for every $i\in[k-1]$,
\(|\Pi\cap H_i|\ge \frac{N-k}{2k}\) and \(|\Pi\cap V_i|\ge \frac{N-k}{2k}\),
where
\(H_i=\RR\times(i,i+1)\) and \(V_i=(i,i+1)\times\RR\).
\end{theorem}

Theorem~\ref{thm:geometric-construction} implies Theorem~1.2. Let $\Pi$ and $\pi$ be supplied by Theorem~\ref{thm:geometric-construction}. After the points of $\Pi$ are replaced by their ranks in the $x$- and $y$-orders, the rank gap between the two canonical points on opposite sides of $V_i$ is at least $|\Pi\cap V_i|+1$, and the analogous rank gap across $H_i$ is at least $|\Pi\cap H_i|+1$. Hence every two distinct points of the unique copy $\Sigma$ are separated in each rank coordinate by at least
\(\frac{N-k}{2k}\).
Set $\eps_\sigma=1/(8k)$. For $N\ge2k-1$,
\(2\eps_\sigma(N-1)=\frac{N-1}{4k}\le \frac{N-k}{2k}\).
Fact~2.4 therefore shows that $\pi$ is $\eps_\sigma$-far from being $\sigma$-free, while Theorem~\ref{thm:geometric-construction} gives exactly one $\sigma$-copy. This proves Theorem~1.2 for every $N\ge\max\{N_0(k),2k-1\}$.

We now prove Theorem~\ref{thm:geometric-construction}. Fix $\sigma\in S_k$ and its canonical point set $\Sigma$. In each case below, we construct finitely many line-segment subsets $P_1,\ldots,P_m$ and set
\(D=\Sigma\cup P_1\cup\cdots\cup P_m\).
The corresponding structural lemma shows that every point set $\Pi$ satisfying
\(\Sigma\subseteq\Pi\subseteq D\)
and associated with a permutation contains $\Sigma$ as its unique $\sigma$-copy.

The segment family in each case also has the following coverage property: every open horizontal stripe $H_i$ and every open vertical stripe $V_i$ contains a nonempty open subsegment of one of the $P_j$. Choose at most one such subsegment for each of the $2(k-1)$ stripes, merge repetitions, and shorten them if necessary so that the chosen subsegments are pairwise disjoint. Distribute the $N-k$ new points as evenly as possible among these subsegments. The points can be chosen successively so that no two have the same $x$- or $y$-coordinate and no new point lies on an integer coordinate line. Thus the resulting set is associated with a permutation and each selected subsegment receives at least
\(\left\lfloor\frac{N-k}{2(k-1)}\right\rfloor\)
points. For $N-k\ge2k(k-1)$,
\(\left\lfloor\frac{N-k}{2(k-1)}\right\rfloor \ge \frac{N-k}{2k}\).
Consequently, once the structural lemmas are proved, the stripe bounds in Theorem~\ref{thm:geometric-construction} follow for all sufficiently large $N$.

For $i,j\in[k]$, write $\Sigma(i)=(i,\sigma(i))$ and $\Sigma^{-1}(j)=(\sigma^{-1}(j),j)$. Because the constructions depend repeatedly on the four boundary points, define $\ell=\sigma(1)$, $r=\sigma(k)$, $b=\sigma^{-1}(1)$, and $u=\sigma^{-1}(k)$, and write $L=\Sigma(1)=(1,\ell)$, $R=\Sigma(k)=(k,r)$, $B=\Sigma^{-1}(1)=(b,1)$, and $U=\Sigma^{-1}(k)=(u,k)$.
Whenever $\Theta$ is a $\sigma$-copy, $\theta$ denotes the permutation associated with $\Theta$; thus $\theta\cong\sigma$.

\begin{definition}\label{def:slope}
A set $P\subset\RR^2$ has \emph{positive slope} if, for every two distinct points $p,q\in P$,
\(p_x<q_x\quad\Leftrightarrow\quad p_y<q_y\).
It has \emph{negative slope} if, for every two distinct points $p,q\in P$,
\(p_x<q_x\quad\Leftrightarrow\quad p_y>q_y\).
\end{definition}

Every subset of a line of direction $(1,1)$ has positive slope, and every subset of a line of direction $(1,-1)$ has negative slope. The definition applies equally to a segment with finitely many points removed and to a union of collinear open subsegments.

The following boundary observation will be used repeatedly.

\begin{claim}\label{clm:boundary}
Let $D\supseteq\Sigma$, and let $\Theta\subseteq D$ be a $k$-point set associated with a permutation isomorphic to $\sigma$.
\begin{enumerate}[label=\textup{(\arabic*)},leftmargin=2.5em]
\item Suppose $P_L=\{p\in D:p_x<1\}$. If $P_L$ has negative slope when $\sigma(1)<\sigma(2)$ and positive slope when $\sigma(1)>\sigma(2)$, then
\(\Theta\cap P_L\subseteq\{\Theta(1)\}\).
\item Suppose $P_R=\{p\in D:p_x>k\}$. If $P_R$ has negative slope when $\sigma(k-1)<\sigma(k)$ and positive slope when $\sigma(k-1)>\sigma(k)$, then
\(\Theta\cap P_R\subseteq\{\Theta(k)\}\).
\item Suppose $P_B=\{p\in D:p_y<1\}$. If $P_B$ has negative slope when $\sigma^{-1}(1)<\sigma^{-1}(2)$ and positive slope when $\sigma^{-1}(1)>\sigma^{-1}(2)$, then
\(\Theta\cap P_B\subseteq\{\Theta^{-1}(1)\}\).
\item Suppose $P_U=\{p\in D:p_y>k\}$. If $P_U$ has negative slope when $\sigma^{-1}(k-1)<\sigma^{-1}(k)$ and positive slope when $\sigma^{-1}(k-1)>\sigma^{-1}(k)$, then
\(\Theta\cap P_U\subseteq\{\Theta^{-1}(k)\}\).
\end{enumerate}
\end{claim}

\begin{proof}
We prove (1); the other statements follow by the corresponding reflections of the plane. Suppose $\Theta(j)\in P_L$ for some $j\ge2$. Since
\(\Theta(1)_x<\Theta(2)_x<\cdots<\Theta(j)_x<1\)
and $P_L$ contains every point of $D$ with $x$-coordinate less than $1$, both $\Theta(1)$ and $\Theta(2)$ lie in $P_L$. If $\sigma(1)<\sigma(2)$, then these two points must increase in both coordinates, contradicting the negative slope of $P_L$. If $\sigma(1)>\sigma(2)$, then they must increase in $x$ and decrease in $y$, contradicting the positive slope of $P_L$. Thus no such $j$ exists.
\end{proof}

\subsection{Patterns with two monotone blocks}\label{subsec:two-blocks}

Assume that $\sigma$ has one of the forms
\((\ell(\ell-1)\cdots21)(k(k-1)\cdots(r+1)r)\), \(r=\ell+1\),
or
\((\ell(\ell+1)\cdots(k-1)k)(12\cdots(r-1)r)\), \(\ell=r+1\).
The two families lie in the same dihedral symmetry orbit. It therefore suffices to treat
\(\sigma=(\ell(\ell-1)\cdots21)(k(k-1)\cdots(r+1)r)\), \(r=\ell+1\).

\subsubsection[Boundary case]{The boundary case $\ell=1$ or $r=k$}

By symmetry, assume $\ell=1$, so
\(\sigma=(1\,k\,(k-1)\cdots3\,2)\), \(k\ge4\).

\begin{construction}\label{con:boundary-two-block}
Let
\(P_0=(0,k+1)+t(-1,-1)\), \(0<t<1\),
and, for $i\in\ZZ$, let
\(P_i=P_0+(i,-i)\).
Set
\(P_M=\bigcup\limits_{i=1}^{k}P_i\).
Thus $P_i$ has coordinate ranges
\(p_x\in(i-1,i)\), \(p_y\in(k-i,k-i+1)\).
\end{construction}

\begin{figure}[H]
\centering
\begin{tikzpicture}[scale=.53]
 \draw[step=1,gray!35,very thin] (0,0) grid (8,8);
 \draw[->] (-.35,0)--(8.45,0);
 \draw[->] (0,-.35)--(0,8.45);
 \foreach \i in {1,...,8}{
 \pgfmathtruncatemacro{\y}{10-\i}
 \ifnum\i=1\relax\pgfmathtruncatemacro{\y}{1}\fi
 \fill[red] (\i,\y) circle (2.2pt);
 }
 \foreach \i in {1,...,8}{
 \pgfmathsetmacro{\xa}{\i-.92}
 \pgfmathsetmacro{\ya}{8-\i+.08}
 \pgfmathsetmacro{\xb}{\i-.08}
 \pgfmathsetmacro{\yb}{8-\i+.92}
 \draw[blue,densely dotted,very thick] (\xa,\ya)--(\xb,\yb);
 }
 \node[below] at (4,-.55) {$x$};
 \node[left] at (-.45,4) {$y$};
\end{tikzpicture}
\caption{Construction~\ref{con:boundary-two-block} for $k=8$. The red points form $\Sigma$ for $\sigma=(1,8,7,6,5,4,3,2)$, and the dotted segments are the sets $P_i$.}
\label{fig:boundary-two-block}
\end{figure}
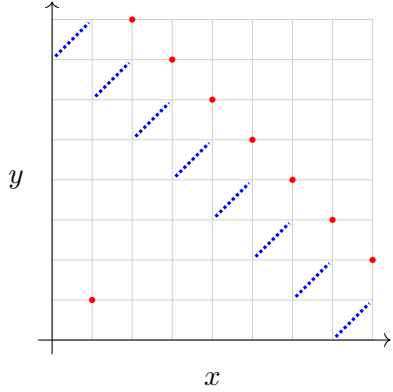

\begin{lemma}\label{lem:boundary-two-block}
If $\Pi$ is associated with a permutation and
\(\Sigma\subseteq\Pi\subseteq\Sigma\cup P_M\),
then $\Sigma$ is the unique $\sigma$-copy in $\Pi$.
\end{lemma}

\begin{proof}
Set $D=\Sigma\cup P_M$, and suppose that $\Theta\subseteq D$ is a $\sigma$-copy distinct from $\Sigma$.

\begin{claim}\label{clm:first-fixed}
$\Theta(1)=\Sigma(1)$.
\end{claim}

\begin{proof}
Assume first that $\Theta(1)=\Sigma(j)$ for some $j\ge2$. No point of $D$ has both coordinates strictly larger than those of $\Sigma(j)$: a canonical point to its right lies below it, while every segment to its right also lies below it. This contradicts the fact that every other point of a copy of $(1\,k\,(k-1)\cdots2)$ lies strictly northeast of its first point.

Now suppose $\Theta(1)\in P_j$. The only points of $D$ lying strictly northeast of $\Theta(1)$ are points of $P_j$ together with at most the two canonical points $\Sigma(j)$ and $\Sigma(j+1)$ when those points exist. Since $k\ge4$, at least one point $\Theta(h)$ with $2\le h\le k$ also lies in $P_j$. Because
\(\theta(2)\ge\theta(h)\ge\theta(k)>\theta(1)\),
the $x$-coordinate inequalities force $\Theta(2)\in P_j$, and the $y$-coordinate inequalities force $\Theta(k)\in P_j$. But $P_j$ has positive slope, whereas $2<k$ and $\theta(2)>\theta(k)$. This contradiction proves the claim.
\end{proof}

Let
\(\Theta_U=\{\Theta(i):2\le i\le k\}\).
For $2\le i<i'\le k$, the pattern $\sigma$ gives
\(\Theta(i)_x<\Theta(i')_x\) and \(\Theta(i)_y>\Theta(i')_y\).
Thus no two points of $\Theta_U$ increase in both coordinates.

Because $\Theta\ne\Sigma$, at least one point of $\Theta_U$ lies on some segment $P_j$.

\begin{claim}\label{clm:end-segments-empty}
$\Theta_U\cap P_1=\Theta_U\cap P_k=\varnothing$.
\end{claim}

\begin{proof}
Since $\Theta(1)=\Sigma(1)=(1,1)$ and every other point of $\Theta$ lies strictly northeast of it, no point of $\Theta_U$ can lie on $P_1$, whose points have $x<1$, or on $P_k$, whose points have $y<1$.
\end{proof}

Choose $j\in\{2,\ldots,k-1\}$ with $\Theta_U\cap P_j\ne\varnothing$. For $2\le s\le k-1$, define $S_s=P_s\cup\{\Sigma(s)\}$ when $s<j$, define $S_j=P_j\cup\{\Sigma(j),\Sigma(j+1)\}$, and define $S_s=P_s\cup\{\Sigma(s+1)\}$ when $s>j$.
Then
\(D=\{\Sigma(1)\}\cup P_1\cup P_k\cup\bigcup\limits_{s=2}^{k-1}S_s\).

\begin{claim}\label{clm:one-per-cell}
For every $2\le s\le k-1$, one has $|\Theta_U\cap S_s|\le1$.
\end{claim}

\begin{proof}
For $s\ne j$, the set $S_s$ has positive slope, so it cannot contain two points of the decreasing set $\Theta_U$. For $s=j$, fix $q\in\Theta_U\cap P_j$. Every other point of $S_j$ is strictly northeast or southwest of $q$ in the positive-slope order and hence cannot also belong to $\Theta_U$. Thus $|\Theta_U\cap S_j|\le1$ as well.
\end{proof}

Claims~\ref{clm:end-segments-empty} and~\ref{clm:one-per-cell} now give
\(|\Theta_U|\le\sum\limits_{s=2}^{k-1}|\Theta_U\cap S_s|\le k-2\),
contradicting $|\Theta_U|=k-1$. Thus $\Sigma$ is the unique copy.
\end{proof}

\subsubsection[Interior case]{The interior case $2\le\ell,r\le k-1$}

By reverse-complement symmetry, assume $\ell\le k/2$.

\begin{construction}\label{con:interior-two-block}
Let
\(P_0=(\ell,\ell)+t(1,1)\), \(0<t<1\),
and let
\(P_i=P_0+(-i,i)\)
for every integer $i$. Set
\(P_M=\bigcup\limits_{i=-k+\ell+1}^{k-\ell-1}P_i\).
Each $P_i$ has coordinate ranges
\(p_x\in(\ell-i,\ell-i+1)\), \(p_y\in(\ell+i,\ell+i+1)\).
\end{construction}

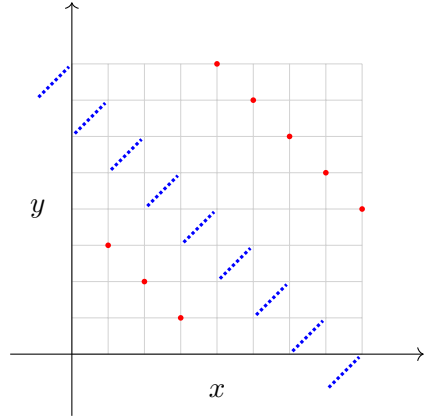
\begin{figure}[H]
\centering
\begin{tikzpicture}[scale=.48]
 \draw[step=1,gray!35,very thin] (0,0) grid (8,8);
 \draw[->] (-1.7,0)--(9.7,0);
 \draw[->] (0,-1.7)--(0,9.7);
 \foreach \i/\y in {1/3,2/2,3/1,4/8,5/7,6/6,7/5,8/4}
 \fill[red] (\i,\y) circle (2.2pt);
 \foreach \j in {-4,...,4}{
 \pgfmathsetmacro{\xa}{3-\j+.08}
 \pgfmathsetmacro{\ya}{3+\j+.08}
 \pgfmathsetmacro{\xb}{3-\j+.92}
 \pgfmathsetmacro{\yb}{3+\j+.92}
 \draw[blue,densely dotted,very thick] (\xa,\ya)--(\xb,\yb);
 }
 \node[below] at (4,-.55) {$x$};
 \node[left] at (-.45,4) {$y$};
\end{tikzpicture}
\caption{Construction~\ref{con:interior-two-block} for $k=8$ and $\ell=3$.}
\label{fig:interior-two-block}
\end{figure}

\begin{lemma}\label{lem:interior-two-block}
If $\Pi$ is associated with a permutation and
\(\Sigma\subseteq\Pi\subseteq\Sigma\cup P_M\),
then $\Sigma$ is the unique $\sigma$-copy in $\Pi$.
\end{lemma}

\begin{proof}
Set $D=\Sigma\cup P_M$, and suppose that $\Theta\subseteq D$ is a $\sigma$-copy distinct from $\Sigma$. Then $\Theta(a)\in P_j$ for some $a\in[k]$ and some
\(-k+\ell+1\le j\le k-\ell-1\).
We first treat $a\le\ell$. Let
\(\Theta_U=\{\Theta(i):r\le i\le k\}\).
Every point of $\Theta_U$ lies strictly northeast of $\Theta(a)$.

\begin{claim}\label{clm:upper-block-canonical}
$\Theta_U\cap P_M=\varnothing$.
\end{claim}

\begin{proof}
Suppose $\Theta(h)\in P_{j'}$ for some $h\ge r$. If $j'\ne j$, then the translated segments $P_j$ and $P_{j'}$ are coordinate-incomparable: no point of either segment lies strictly northeast of a point of the other. This contradicts $\Theta(h)$ lying northeast of $\Theta(a)$.

If $j'=j$, then $\Theta(a),\Theta(h)\in P_j$. Since
\(a<r\le h\le k\) and \(\theta(a)<\theta(k)\le\theta(h)\le\theta(r)\),
the $x$-coordinate range of $P_j$ forces $\Theta(r)\in P_j$, while its $y$-coordinate range forces $\Theta(k)\in P_j$. This is impossible because $P_j$ has positive slope but $r<k$ and $\theta(r)>\theta(k)$.
\end{proof}

Thus $\Theta_U\subseteq\Sigma$. Since $|\Theta_U|=k-\ell$ and every point of $\Theta_U$ lies strictly northeast of $\Theta(a)$, there must be at least $k-\ell$ canonical $x$-coordinates and at least $k-\ell$ canonical $y$-coordinates larger than the corresponding coordinates of $\Theta(a)$. Hence
\(\Theta(a)_x<\ell+1\) and \(\Theta(a)_y<\ell+1\).
Because $\Theta(a)\in P_j$, these inequalities give $-1<j<1$, so $j=0$. Therefore $\Theta(a)\in P_0\subset(\ell,\ell+1)^2$. The only $k-\ell$ points of $\Sigma$ lying strictly northeast of this square are $\Sigma(r),\ldots,\Sigma(k)$, so
\(\Theta(i)=\Sigma(i)\), \((r\le i\le k)\).

Now
\(\ell<\Theta(a)_x\le\Theta(\ell)_x<\Theta(r)_x=r=\ell+1\)
and
\(\ell<\Theta(a)_y\le\Theta(1)_y<\Theta(k)_y=r=\ell+1\).
The only points of $D$ with an $x$- or $y$-coordinate in $(\ell,\ell+1)$ lie on $P_0$. Hence $\Theta(\ell),\Theta(1)\in P_0$, contradicting the positive slope of $P_0$ because $1<\ell$ and $\theta(1)>\theta(\ell)$.

If $a\ge r$, the same argument with southwest in place of northeast first shows that the lower block $\{\Theta(i):1\le i\le\ell\}$ lies in $\Sigma$, then forces $j=0$, and finally places $\Theta(r)$ and $\Theta(k)$ together on $P_0$, again contradicting its positive slope. Thus no such $\Theta$ exists.
\end{proof}

\subsection[The patterns 31524 and 24153]{The patterns $(31524)$, $(24153)$, $(42513)$, and $(35142)$}\label{subsec:exceptional-five}

The four patterns are related by a dihedral symmetry, so it suffices to take $\sigma=(31524)$.

\begin{construction}\label{con:31524}
Let
\(I_M=(1,5)+t(1,-1),\quad I_B=(3,1)+t(1,-1),\quad I_R=(5,3)+t(1,-1)\), \(0\le t\le2\).
Define
\(P_M=I_M\setminus\ZZ^2\), \(P_B=I_B\setminus\ZZ^2\), \(P_R=I_R\setminus\ZZ^2\).
\end{construction}

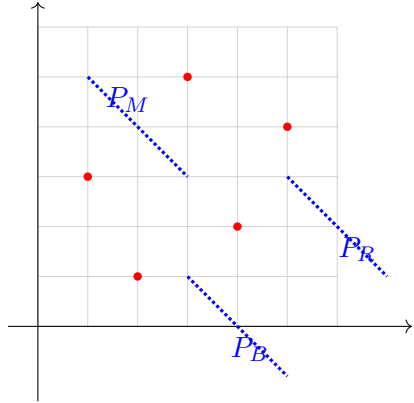
\begin{figure}[H]
\centering
\begin{tikzpicture}[scale=.66]
 \draw[step=1,gray!35,very thin] (0,0) grid (6,6);
 \draw[->] (-.6,0)--(7.5,0);
 \draw[->] (0,-1.5)--(0,6.5);
 \foreach \i/\y in {1/3,2/1,3/5,4/2,5/4}
 \fill[red] (\i,\y) circle (2.4pt);
 \draw[blue,densely dotted,very thick] (1,5)--(3,3);
 \draw[blue,densely dotted,very thick] (3,1)--(5,-1);
 \draw[blue,densely dotted,very thick] (5,3)--(7,1);
 \node[blue] at (1.8,4.55) {$P_M$};
 \node[blue] at (4.25,-.45) {$P_B$};
 \node[blue] at (6.4,1.55) {$P_R$};
\end{tikzpicture}
\caption{Construction~\ref{con:31524}. Integer lattice points on the three dotted lines are deleted.}
\label{fig:31524}
\end{figure}

\begin{lemma}\label{lem:31524}
If $\Pi$ is associated with a permutation and
\(\Sigma\subseteq\Pi\subseteq\Sigma\cup P_M\cup P_B\cup P_R\),
then $\Sigma$ is the unique $\sigma$-copy in $\Pi$.
\end{lemma}

\begin{proof}
Set $D=\Sigma\cup P_M\cup P_B\cup P_R$, and suppose that $\Theta\subseteq D$ is a $\sigma$-copy distinct from $\Sigma$. The sets $P_R$ and $P_B$ satisfy parts (2) and (3), respectively, of Claim~\ref{clm:boundary}. Therefore
\(\Theta\cap P_R\subseteq\{\Theta(5)\}\), \(\Theta\cap P_B\subseteq\{\Theta^{-1}(1)\}=\{\Theta(2)\}\).

\begin{claim}\label{clm:no-bottom-right}
$\Theta\cap P_B=\Theta\cap P_R=\varnothing$.
\end{claim}

\begin{proof}
Assume $\Theta\cap P_R\ne\varnothing$. Then $\Theta(5)\in P_R$, so $\Theta(5)_y<3$. Since $\theta(1)<\theta(5)$, also $\Theta(1)_y<3$. All points of $D$ with $y<3$ lie in
\(P_B\cup P_R\cup\{\Sigma(2),\Sigma(4)\}\).
Claim~\ref{clm:boundary} excludes $P_B$ and $P_R$ for $\Theta(1)$, so $\Theta(1)\in\{\Sigma(2),\Sigma(4)\}$. Because $\theta(4)<\theta(1)$ and $4>1$, the point $\Theta(4)$ lies strictly southeast of $\Theta(1)$. Every point of $D$ satisfying the resulting coordinate inequalities lies in $P_B\cup P_R$, but Claim~\ref{clm:boundary} allows neither boundary set to contain $\Theta(4)$. This contradiction proves $\Theta\cap P_R=\varnothing$.

Reflection in the anti-diagonal $x+y=6$ preserves $D$, fixes the pattern $(31524)$, and interchanges $P_B$ with $P_R$. Applying the argument to the reflected copy gives $\Theta\cap P_B=\varnothing$.
\end{proof}

It follows that $\Theta\subseteq\Sigma\cup P_M$ and $\Theta\cap P_M\ne\varnothing$.

If $\Theta(4)\in P_M$, then $\Theta(4)_x<3$ and $\Theta(4)_y>3$. Since
\(\theta(3)>\theta(1)>\theta(4)\),
both $\Theta(1)$ and $\Theta(3)$ also have $x<3$ and $y>3$, so both lie in $P_M$. But $\Theta(1)_x<\Theta(3)_x$ and $\Theta(1)_y<\Theta(3)_y$, contradicting the negative slope of $P_M$.

Otherwise, reflecting the copy in the anti-diagonal if necessary allows us to assume that $\Theta(1)\in P_M$ or $\Theta(2)\in P_M$. In either case, both $\Theta(3)$ and $\Theta(5)$ lie strictly northeast of a point of $P_M$. They cannot themselves lie in $P_M$, because a negative-slope set contains no two points that increase in both coordinates.

Thus $\Theta(3),\Theta(5)\in\Sigma$. They are the two points of $\Sigma$ with $y>3$, and their $y$-order gives
\(\Theta(3)=\Sigma(3)\), \(\Theta(5)=\Sigma(5)\).
It follows successively that $\Theta(4)=\Sigma(4)$, because it is the only available point with $3<x<5$; that $\Theta(2)=\Sigma(2)$, because it is the only available point with $y<2$; and that $\Theta(1)=\Sigma(1)$, because it is the only available point with $x<2$ and $y<4$. Hence $\Theta=\Sigma$, a contradiction.
\end{proof}

\subsection{Patterns with a coincident boundary point}\label{subsec:coincident-boundary}

Assume
\(\{\sigma(1),\sigma(k)\}\cap\{1,k\}\ne\varnothing\),
equivalently $|\{L,R,B,U\}|\le3$, and assume that $\sigma$ is not covered by Subsections~\ref{subsec:two-blocks} or~\ref{subsec:exceptional-five}. By a dihedral symmetry, we may assume
\(L=B\), \(\text{that is,}\), \(\sigma(1)=1\).

\begin{construction}\label{con:coincident-boundary}
Define four line segments as follows. If $\sigma(2)\ne k$, let $I_L=(1,2)+t(-1,1)$ for $0\le t\le k-2$ and $I_U=(2,k)+t(-1,1)$ for $0\le t\le1$; if $\sigma(2)=k$, let $I_L=(1,3)+t(-1,1)$ for $0\le t\le k-3$ and $I_U=(1,2)+t(1,1)$ for $0\le t\le1$. Similarly, if $\sigma(k)\ne2$, let $I_B=(2,1)+t(1,-1)$ for $0\le t\le k-2$ and $I_R=(k,2)+t(1,-1)$ for $0\le t\le1$; if $\sigma(k)=2$, let $I_B=(3,1)+t(1,-1)$ for $0\le t\le k-3$ and $I_R=(2,1)+t(1,1)$ for $0\le t\le1$. For $X\in\{L,U,B,R\}$, set $P_X=I_X\setminus\ZZ^2$.
\end{construction}

\begin{figure}[H]
\centering
\begin{tikzpicture}[scale=.48]
 \draw[step=1,gray!35,very thin] (1,1) grid (8,8);
 \draw[->] (-5.7,0)--(9.7,0);
 \draw[->] (0,-5.7)--(0,9.7);
 \foreach \i/\y in {1/1,2/4,3/7,4/3,5/8,6/5,7/2,8/6}
 \fill[red] (\i,\y) circle (2.3pt);
 \draw[blue,densely dotted,very thick] (1,2)--(-5,8);
 \draw[blue,densely dotted,very thick] (2,8)--(1,9);
 \draw[blue,densely dotted,very thick] (2,1)--(8,-5);
 \draw[blue,densely dotted,very thick] (8,2)--(9,1);
 \node[blue] at (-2.7,6.2) {$P_L$};
 \node[blue] at (1.15,9.25) {$P_U$};
 \node[blue] at (6,-3.8) {$P_B$};
 \node[blue] at (9.2,1.45) {$P_R$};
\end{tikzpicture}
\caption{Construction~\ref{con:coincident-boundary} in the branches $\sigma(2)\ne k$ and $\sigma(k)\ne2$. In the exceptional branches, the short segment $P_U$ or $P_R$ is moved into the indicated missing unit stripe and has positive slope.}
\label{fig:coincident-boundary}
\end{figure}
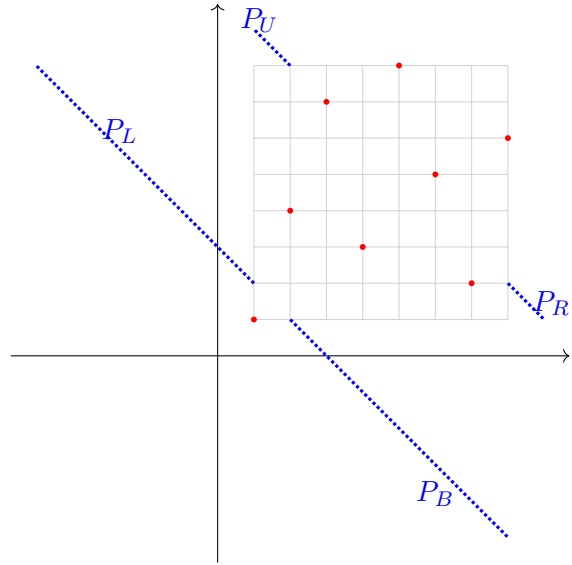

\begin{lemma}\label{lem:coincident-boundary}
If $\Pi$ is associated with a permutation and
\(\Sigma\subseteq\Pi\subseteq\Sigma\cup P_L\cup P_R\cup P_B\cup P_U\),
then $\Sigma$ is the unique $\sigma$-copy in $\Pi$.
\end{lemma}

\begin{proof}
Set $D=\Sigma\cup P_L\cup P_R\cup P_B\cup P_U$, and suppose that $\Theta\subseteq D$ is a $\sigma$-copy. Since $\sigma(1)=1$, we have
\(\Theta(1)=\Theta^{-1}(1)\).
The sets $P_L$ and $P_B$ satisfy parts (1) and (3) of Claim~\ref{clm:boundary}, respectively, so
\(\Theta\cap P_L\subseteq\{\Theta(1)\}\), \(\Theta\cap P_B\subseteq\{\Theta(1)\}\).

\begin{claim}\label{clm:canonical-corner-exclusion}
If $\Theta\cap P_L\ne\varnothing$, then $\Sigma(1)\notin\Theta$. Likewise, if $\Theta\cap P_B\ne\varnothing$, then $\Sigma(1)\notin\Theta$.
\end{claim}

\begin{proof}
If $\Theta(1)\in P_L$ and $\Sigma(1)=\Theta(j)$ for some $j>1$, then
\(\Theta(1)_y>1=\Theta(j)_y\),
contradicting $\theta(1)=1<\theta(j)$. If $\Theta(1)\in P_B$ and $\Sigma(1)=\Theta(j)$ for some $j>1$, then
\(\Theta(1)_x>1=\Theta(j)_x\),
contradicting the $x$-order $\Theta(1)_x<\Theta(j)_x$.
\end{proof}

\begin{claim}\label{clm:opposite-segments}
The set $\Theta$ cannot meet both $P_L$ and $P_U$. Likewise, it cannot meet both $P_B$ and $P_R$.
\end{claim}

\begin{proof}
Assume that $\Theta(1)\in P_L$ and $\Theta\cap P_U\ne\varnothing$. By Claim~\ref{clm:canonical-corner-exclusion}, $\Sigma(1)\notin\Theta$. All points of $D$ other than $\Sigma(1)$, points of $P_L$ (which contains $\Theta(1)$) and points of $P_U$ have $x$-coordinate at least $2$, whereas $P_U$ has $x$-coordinate in $(1,2)$. Hence $\Theta(2)$ is the leftmost point of $\Theta\cap P_U$.

If $\sigma(2)=k$, then $P_L$ has $y>3$ while $P_U$ has $2<y<3$, so $\Theta(1)_y>\Theta(2)_y$, contradicting $\theta(1)<\theta(2)$. If $\sigma(2)\ne k$, then $P_U$ has $y>k$ and negative slope. Its leftmost selected point has the largest $y$-coordinate among $\Theta\cap P_U$, and every point outside $P_U$ has $y\le k$. Thus $\Theta(2)=\Theta^{-1}(k)$, forcing $\sigma(2)=k$, again a contradiction.

Now assume $\Theta(1)\in P_B$ and $\Theta\cap P_R\ne\varnothing$. The coordinate-dual argument identifies $\Theta^{-1}(2)$ as the lowest point of $\Theta\cap P_R$. If $\sigma(k)=2$, then $P_B$ has $x>3$ while $P_R$ has $2<x<3$, contradicting the fact that $\Theta(1)$ precedes $\Theta(k)=\Theta^{-1}(2)$. If $\sigma(k)\ne2$, then $P_R$ has $x>k$ and negative slope; its lowest selected point is the rightmost point of $\Theta$, so $\Theta^{-1}(2)=\Theta(k)$, forcing $\sigma(k)=2$. This is again a contradiction.
\end{proof}

We now show that $\Theta$ cannot meet $P_L\cup P_U$. First suppose $\Theta\cap P_L\ne\varnothing$. Then $\Theta(1)\in P_L$, Claim~\ref{clm:opposite-segments} gives $\Theta\cap P_U=\varnothing$, and $\Theta(1)_y>2$. Since $\Theta(1)$ is the lowest point of $\Theta$, no point of $\Theta$ lies in $P_B\cup P_R$, whose points have $y<2$. Moreover, only $k-2$ canonical points have $y>2$. Hence
\(|\Theta|\le |\Theta\cap P_L|+|\Theta\cap\Sigma| \le1+(k-2)=k-1\),
a contradiction.

It remains to consider $\Theta\cap P_L=\varnothing$ and $\Theta\cap P_U\ne\varnothing$. We first prove $\Theta(1)=\Sigma(1)$. If not, then $\Sigma(1)\notin\Theta$: because no point of $P_L$ is selected, $\Sigma(1)$ would be the leftmost point of $\Theta$ if it were selected. Since $P_U$ is nonempty in $\Theta$ and has $1<x<2$, while every selected point outside $P_U$ has $x\ge2$, the leftmost point $\Theta(1)$ lies in $P_U$.

Note that $k\ge4$, because the only length-three pattern with $\sigma(1)=1$, $\sigma(2)=3$, and $\sigma(3)=2$ is $(132)$, already covered by Subsection~\ref{subsec:two-blocks}. 

If $\sigma(2)\ne k$, then $P_U$ has negative slope and lies above $y=k$. No point to the right of $\Theta(1)$ can lie above it: a point of $P_U$ to the right lies below it, and every point outside $P_U$ has $y\le k$. This contradicts the fact that $\Theta(1)$ is the minimum point in the $y$-order.

If $\sigma(2)=k$, then $P_U\subset(1,2)\times(2,3)$ has positive slope. Since $\Theta(1)_y>2$, all $k$ points of $\Theta$ must have $y>2$. Outside $P_U$ and $P_L$, there are only the $k-2$ canonical points of values $3,\ldots,k$. Thus at least two points of $\Theta$ lie in $P_U$, and they are the first two points in the $x$-order. If a canonical point of value at least $3$ is used, it lies above every point of $P_U$, contradicting that $\Theta(2)$ is the maximum point of the copy. If no canonical point is used, all points lie on the positive-slope set $P_U$, so a later point lies above $\Theta(2)$, giving the same contradiction. Therefore $\Theta(1)=\Sigma(1)$.

Because $P_U$ is nonempty in $\Theta$ and every point outside $P_U$ other than $\Sigma(1)$ has $x\ge2$, the point $\Theta(2)$ lies in $P_U$. If $\sigma(2)\ne k$, the same leftmost-on-a-negative-slope argument makes $\Theta(2)$ the maximum point of $\Theta$, a contradiction. Hence
\(\sigma(2)=k\), \(P_U\subset(1,2)\times(2,3)\), \(\Theta(2)=\Theta^{-1}(k)\).
No further point of $\Theta$ can lie in $P_U$: a point to the right on this positive-slope segment would lie above $\Theta(2)$.

All remaining points $\Theta(3),\ldots,\Theta(k)$ have $1<y<3$. The only available points in this strip are $\Sigma^{-1}(2)$ and points of $P_R$. If $\sigma(k)=2$, then $P_R$ has positive slope and lies in $(2,3)\times(1,2)$, while $\Sigma^{-1}(2)=\Sigma(k)$ lies northeast of every point of $P_R$. Thus the entire available set $P_R\cup\{\Sigma^{-1}(2)\}$ has positive slope. Since $k\geq4$, the tail $\Theta(3),\ldots,\Theta(k)$ therefore contains at least two points and would increase in $y$, contradicting that $\theta(k)=2$ is its smallest value.

Thus $\sigma(k)\ne2$. In this branch $P_R$ lies to the right of $x=k$ and has negative slope, while $\Sigma^{-1}(2)$ lies northwest of every point of $P_R$. Hence the entire available set $P_R\cup\{\Sigma^{-1}(2)\}$ has negative slope, and therefore
\(\theta(3)>\theta(4)>\cdots>\theta(k)\),
so
\(\sigma=(1\,k\,(k-1)\cdots3\,2)\),
which is excluded by Subsection~\ref{subsec:two-blocks} (for $k\geq4$) and by the hypothesis of Theorem~\ref{thm:geometric-construction} (for $k=3$). This contradiction proves
\(\Theta\cap P_L=\Theta\cap P_U=\varnothing\).

Reflecting the entire construction across the line $y=x$ replaces $\sigma$ by $\sigma^{-1}$, interchanges $P_L$ with $P_B$ and $P_U$ with $P_R$, and preserves the standing exclusions because the earlier pattern families are closed under inversion. Applying the argument above to the reflected copy gives
\(\Theta\cap P_B=\Theta\cap P_R=\varnothing\).
Therefore $\Theta\subseteq\Sigma$. Since both sets have $k$ points, $\Theta=\Sigma$. This proves the lemma.
\end{proof}

\subsection{Patterns where the four boundary points are distinct}\label{subsec:four-boundary-points}

Throughout this subsection, assume that the pattern $\sigma\in S_k$ is not covered by Subsections~\ref{subsec:two-blocks} or~\ref{subsec:exceptional-five} and that the four boundary points
\(L=\Sigma(1)\), \(R=\Sigma(k)\), \(B=\Sigma^{-1}(1)\), \(U=\Sigma^{-1}(k)\)
are pairwise distinct. Equivalently,
\(\{\sigma(1),\sigma(k)\}\cap\{1,k\}=\varnothing\),
so
\(1<\ell,r,b,u<k\) with $k\geq4$, \(\ell=\sigma(1),\ r=\sigma(k),\ b=\sigma^{-1}(1),\ u=\sigma^{-1}(k)\).

\begin{construction}\label{con:four-boundary}
Define four closed diagonal segments according to the local boundary behavior of $\sigma$. Let $I_{\mathrm L}=(1,1)+t(-1,1)$ when $\sigma(2)>\ell$ and $I_{\mathrm L}=(1,k)+t(-1,-1)$ when $\sigma(2)<\ell$, in both cases for $0\le t\le k-1$. Let $I_{\mathrm R}=(k,1)+t(1,1)$ when $\sigma(k-1)>r$ and $I_{\mathrm R}=(k,k)+t(1,-1)$ when $\sigma(k-1)<r$, again for $0\le t\le k-1$. Let $I_{\mathrm B}=(1,1)+t(1,-1)$ when $\sigma^{-1}(2)>b$ and $I_{\mathrm B}=(k,1)+t(-1,-1)$ when $\sigma^{-1}(2)<b$, for $0\le t\le k-1$. Finally, let $I_{\mathrm U}=(1,k)+t(1,1)$ when $\sigma^{-1}(k-1)>u$ and $I_{\mathrm U}=(k,k)+t(-1,1)$ when $\sigma^{-1}(k-1)<u$, for $0\le t\le k-1$.
First delete every lattice point from these four segments. Then trim the remaining sets as follows. If $\ell<r$, let
\(\PL=\{p\in I_{\mathrm L}\setminus\mathbb Z^2:p_y>\ell+1\}\), \(\PR=\{p\in I_{\mathrm R}\setminus\mathbb Z^2:p_y<r-1\}\);
if $\ell>r$, let
\(\PL=\{p\in I_{\mathrm L}\setminus\mathbb Z^2:p_y<\ell-1\}\), \(\PR=\{p\in I_{\mathrm R}\setminus\mathbb Z^2:p_y>r+1\}\).
Similarly, if $b<u$, let
\(\PB=\{p\in I_{\mathrm B}\setminus\mathbb Z^2:p_x>b+1\}\), \(\PU=\{p\in I_{\mathrm U}\setminus\mathbb Z^2:p_x<u-1\}\);
if $b>u$, let
\(\PB=\{p\in I_{\mathrm B}\setminus\mathbb Z^2:p_x<b-1\}\), \(\PU=\{p\in I_{\mathrm U}\setminus\mathbb Z^2:p_x>u+1\}\).
\end{construction}

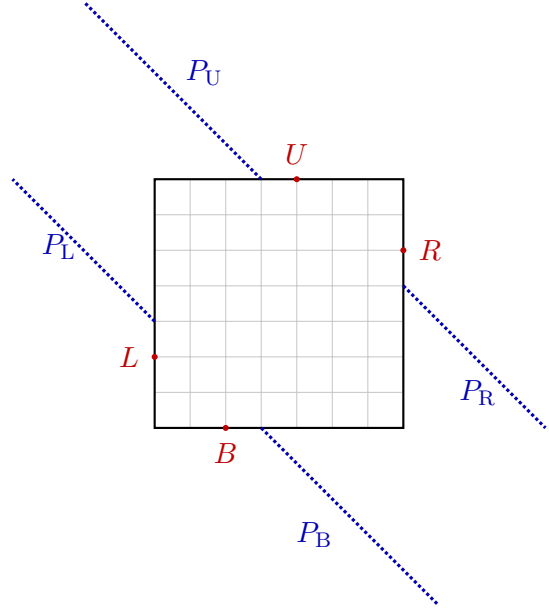
\begin{figure}[H]
\centering
\begin{tikzpicture}[scale=.47]
 \def\K{8}
 \draw[step=1,gray!42,thin] (1,1) grid (8,8);
 \draw[thick] (1,1) rectangle (8,8);
 \draw[blue!70!black,densely dotted,very thick] (1,4)--(-3,8);
 \draw[blue!70!black,densely dotted,very thick] (8,5)--(12,1);
 \draw[blue!70!black,densely dotted,very thick] (4,1)--(9,-4);
 \draw[blue!70!black,densely dotted,very thick] (4,8)--(-1,13);
 \fill[red!75!black] (1,3) circle (2.4pt) node[left=2pt] {$L$};
 \fill[red!75!black] (8,6) circle (2.4pt) node[right=2pt] {$R$};
 \fill[red!75!black] (3,1) circle (2.4pt) node[below=2pt] {$B$};
 \fill[red!75!black] (5,8) circle (2.4pt) node[above=2pt] {$U$};
 \node[blue!70!black] at (-1.7,6.1) {$P_{\mathrm L}$};
 \node[blue!70!black] at (10.1,2.0) {$P_{\mathrm R}$};
 \node[blue!70!black] at (5.5,-2.0) {$P_{\mathrm B}$};
 \node[blue!70!black] at (2.4,11.0) {$P_{\mathrm U}$};
\end{tikzpicture}
\caption{A schematic of Construction~\ref{con:four-boundary}. The exact slopes of the four boundary sets are determined by the adjacent entries of $\sigma$.}
\end{figure}

\begin{lemma}\label{lem:four-boundary-unique}
If $\Pi$ is associated with a permutation and
\(\Sigma\subseteq\Pi\subseteq \Sigma\cup\PL\cup\PR\cup\PB\cup\PU\),
then $\Sigma$ is the unique $\sigma$-copy in $\Pi$.
\end{lemma}

\begin{proof}
Suppose, for contradiction, that $\Theta\subseteq\Pi$ is a $\sigma$-copy with $\Theta\ne\Sigma$. By Construction~\ref{con:four-boundary}, the four boundary sets satisfy the hypotheses of Claim~\ref{clm:boundary}. Hence
\(\Theta\cap\PL\subseteq\{\Theta(1)\}\), \(\Theta\cap\PR\subseteq\{\Theta(k)\}\),
\inlinetag{\Theta\cap\PB\subseteq\{\Theta^{-1}(1)\},\quad\allowbreak \Theta\cap\PU\subseteq\{\Theta^{-1}(k)\}}{4.4}{eq:boundary-restrictions}.
We first record a propagation property for deviations from the canonical coordinate levels.

\begin{claim}\label{clm:propagation}
The following statements hold.
\begin{enumerate}[label=\arabic*.,leftmargin=2.2em]
\item If $j<k$ and $\Theta^{-1}(j)_y>j$, then $\Theta^{-1}(j+1)_y>j+1$.
\item If $j>1$ and $\Theta^{-1}(j)_y<j$, then $\Theta^{-1}(j-1)_y<j-1$.
\item If $i<k$ and $\Theta(i)_x>i$, then $\Theta(i+1)_x>i+1$.
\item If $i>1$ and $\Theta(i)_x<i$, then $\Theta(i-1)_x<i-1$.
\end{enumerate}
\end{claim}

\begin{proof}
We prove statement~1; the other statements follow by reflecting the construction horizontally or vertically, or by interchanging the two coordinates.

Assume $j<k$ and $\Theta^{-1}(j)_y>j$. If $\Theta^{-1}(j)\in\Sigma$, then its $y$-coordinate is an integer, so
\(\Theta^{-1}(j)_y\ge j+1\),
and therefore
\(\Theta^{-1}(j+1)_y>\Theta^{-1}(j)_y\ge j+1\).

Now suppose $\Theta^{-1}(j)\notin\Sigma$. It cannot lie in $\PB$, because its $y$-coordinate exceeds $j\ge1$, and it cannot lie in $\PU$, because $j<k$ and \eqref{eq:boundary-restrictions} allows only $\Theta^{-1}(k)$ in $\PU$. Thus it lies in $\PL$ or $\PR$.

If $\Theta^{-1}(j)\in\PL$, then \eqref{eq:boundary-restrictions} gives $\Theta^{-1}(j)=\Theta(1)$, and hence $j=\ell$. When $\ell<r$, the trimming of $\PL$ gives $\Theta^{-1}(j)_y>\ell+1=j+1$; when $\ell>r$, every point of $\PL$ has $y<\ell-1$, contradicting the hypothesis $\Theta^{-1}(j)_y>j=\ell$.

If $\Theta^{-1}(j)\in\PR$, then \eqref{eq:boundary-restrictions} gives $\Theta^{-1}(j)=\Theta(k)$, and hence $j=r$. When $\ell<r$, every point of $\PR$ has $y<r-1$, contradicting $\Theta^{-1}(j)_y>j=r$; when $\ell>r$, the trimming of $\PR$ gives $\Theta^{-1}(j)_y>r+1=j+1$.

Thus, in every viable case,
\(\Theta^{-1}(j+1)_y>\Theta^{-1}(j)_y\geq j+1\).
This proves statement~1.
\end{proof}

\begin{corollary}\label{cor:propagation}
The following statements hold.
\begin{enumerate}[label=\arabic*.,leftmargin=2.2em]
\item If $\Theta^{-1}(j)_y>j$ for some $j$, then
\(\Theta^{-1}(s)_y>s\), \((j\le s\le k)\).
\item If $\Theta^{-1}(j)_y<j$ for some $j$, then
\(\Theta^{-1}(s)_y<s\), \((1\le s\le j)\).
\item If $\Theta(i)_x>i$ for some $i$, then
\(\Theta(s)_x>s\), \((i\le s\le k)\).
\item If $\Theta(i)_x<i$ for some $i$, then
\(\Theta(s)_x<s\), \((1\le s\le i)\).
\end{enumerate}
\end{corollary}

\begin{proof}
For statement~1, repeatedly apply statement~1 of Claim~\ref{clm:propagation}: once $\Theta^{-1}(s)_y>s$, the claim gives $\Theta^{-1}(s+1)_y>s+1$. The other three statements follow in the same way from the corresponding parts of Claim~\ref{clm:propagation}.
\end{proof}

Because $\Theta\ne\Sigma$, at least one point of $\Theta$ lies in a boundary set. The construction is invariant under the dihedral symmetries of the square, so we may assume
\inlinetag{\Theta\cap\PL\ne\varnothing \quad\allowbreak\text{and}\quad\allowbreak \ell<r}{4.5}{eq:wlog-left}.
Then \eqref{eq:boundary-restrictions} gives $\Theta\cap\PL=\{\Theta(1)\}$.

\begin{claim}\label{clm:left-point-high}
We have
\(\Theta^{-1}(\ell)=\Theta(1)\) and \(\Theta^{-1}(\ell)_y>\ell+1\).
\end{claim}

\begin{proof}
Since $\Theta$ is order-isomorphic to $\sigma$, the point in position $1$ has value-rank $\ell$, so $\Theta(1)=\Theta^{-1}(\ell)$. Under the assumption $\ell<r$, Construction~\ref{con:four-boundary} retains only points of $\PL$ with $y>\ell+1$.
\end{proof}

\begin{claim}\label{clm:right-canonical}
We have
\(\Theta(k)_y>r\) and \(\Theta(k)\in\Sigma\).
\end{claim}

\begin{proof}
Claim~\ref{clm:left-point-high} and Corollary~\ref{cor:propagation} imply
\(\Theta^{-1}(r)_y>r\).
Because $\Theta^{-1}(r)=\Theta(k)$, this gives the first assertion.

The point $\Theta(k)$ cannot lie in $\PR$, because under $\ell<r$ every point of $\PR$ has $y<r-1$. It cannot lie in $\PL$, $\PB$, or $\PU$ by \eqref{eq:boundary-restrictions}, since $k\ne1$, $r\ne1$, and $r\ne k$. Hence $\Theta(k)\in\Sigma$.
\end{proof}

\begin{claim}\label{clm:top-and-bottom}
We have
\(\Theta^{-1}(k)\in\PU\), \(b<u\), \(\Theta^{-1}(k)_x<u-1\),
and
\(\Theta^{-1}(1)_x<b\), \(\Theta^{-1}(1)\in\Sigma\).
\end{claim}

\begin{proof}
By Claim~\ref{clm:left-point-high} and Corollary~\ref{cor:propagation},
\(\Theta^{-1}(k)_y>k\).
The only part of the support with $y>k$ is $\PU$, so $\Theta^{-1}(k)\in\PU$.

Suppose $b>u$. Then Construction~\ref{con:four-boundary} gives $\Theta^{-1}(k)_x>u+1>u$. Since $\Theta^{-1}(k)=\Theta(u)$, statement~3 of Corollary~\ref{cor:propagation} would imply $\Theta(k)_x>k$, contradicting Claim~\ref{clm:right-canonical}. Hence $b<u$.

For $b<u$, every point of $\PU$ has $x<u-1$, so $\Theta^{-1}(k)_x<u-1$. Since $\Theta^{-1}(k)=\Theta(u)$ and $b<u$, statement~4 of Corollary~\ref{cor:propagation} yields
\(\Theta(b)_x<b\).
But $\Theta(b)=\Theta^{-1}(1)$, so $\Theta^{-1}(1)_x<b$.

The point $\Theta^{-1}(1)$ cannot lie in $\PB$, because under $b<u$ every point of $\PB$ has $x>b+1$. It is excluded from $\PL$, $\PR$, and $\PU$ by \eqref{eq:boundary-restrictions}, because its position is $b\notin\{1,k,u\}$. Thus $\Theta^{-1}(1)\in\Sigma$.
\end{proof}

\begin{claim}\label{clm:shift-structure}
The following two statements hold.
\begin{enumerate}[label=\arabic*.,leftmargin=2.2em]
\item We have
\(u-2<\Theta^{-1}(k)_x=\Theta(u)_x<u-1\), \(b=u-1\),
and
\inlinetag{\Theta(i)=\Sigma(i-1) \quad\allowbreak\text{for every }i\in[k]\setminus\{1,u\}}{4.6}{eq:x-shift}.
\item We have
\(\ell+1<\Theta(1)_y=\Theta^{-1}(\ell)_y<\ell+2\), \(r=\ell+1\),
and
\inlinetag{\Theta^{-1}(j)=\Sigma^{-1}(j+1) \quad\allowbreak\text{for every }j\in[k]\setminus\{\ell,k\}}{4.7}{eq:y-shift}.
\end{enumerate}
\end{claim}

\begin{proof}
We first prove statement~1. By Claim~\ref{clm:top-and-bottom}, $\Theta(u)_x<u-1$. Every point $\Theta(i)$ with $1<i<u$ lies in $\Sigma$: the boundary restriction \eqref{eq:boundary-restrictions} permits a noncanonical point only at positions $1,k,b,u$, while Claims~\ref{clm:right-canonical} and~\ref{clm:top-and-bottom} show that the points at positions $k$ and $b$ are canonical. Thus the $u-2$ points
\(\Theta(2),\Theta(3),\ldots,\Theta(u-1)\)
have distinct positive integral $x$-coordinates less than $\Theta(u)_x<u-1$. Consequently their $x$-coordinates are exactly $1,2,\ldots,u-2$, and
\inlinetag{\Theta(i)_x=i-1\quad\allowbreak(2\le i\le u-1)}{4.8}{eq:first-x-block}.
We next identify the canonical points omitted from $\Theta$. Claim~\ref{clm:right-canonical} gives $\Theta(k)\in\Sigma$ and $\Theta(k)_y>r=\Sigma(k)_y$, so $\Theta(k)\ne\Sigma(k)$. Since $\Theta(k)$ is the rightmost point of $\Theta$, the point $\Sigma(k)$ cannot occur elsewhere in $\Theta$. Similarly, Claim~\ref{clm:top-and-bottom} gives $\Theta^{-1}(1)\in\Sigma$ and $\Theta^{-1}(1)_x<b=\Sigma(b)_x$, so $\Theta^{-1}(1)\ne\Sigma(b)$; because $\Theta^{-1}(1)$ is the lowest point of $\Theta$, the point $\Sigma(b)$ cannot occur elsewhere. The only noncanonical points are therefore $\Theta(1)$ and $\Theta(u)$, and hence
\inlinetag{\Theta\cap\Sigma=\Sigma\setminus\{\Sigma(b),\Sigma(k)\}}{4.8a}{eq:canonical-points-used}.
The coordinates in \eqref{eq:first-x-block} include every integer from $1$ to $u-2$. In view of \eqref{eq:canonical-points-used}, the only missing canonical $x$-coordinates are $b$ and $k$. Since $b<u$, this is possible only if $b=u-1$. The available canonical $x$-coordinates are then
\(1,2,\ldots,u-2,u,u+1,\ldots,k-1\).
Matching them in increasing order to the canonical positions of $\Theta$ gives \eqref{eq:x-shift}. Finally, $\Theta(u-1)_x=u-2<\Theta(u)_x<u-1$, proving the remaining interval in statement~1.

For statement~2, \eqref{eq:canonical-points-used} shows that the canonical $y$-coordinates available to $\Theta\cap\Sigma$ are
\inlinetag{[k]\setminus\{1,r\}}{4.9}{eq:available-y}.
The points $\Theta^{-1}(1),\ldots,\Theta^{-1}(\ell-1)$ are canonical. Since all available canonical points occur in $\Theta$, these are precisely the $\ell-1$ smallest levels in \eqref{eq:available-y}, namely $2,3,\ldots,\ell$. Claim~\ref{clm:left-point-high} gives $\Theta^{-1}(\ell)_y>\ell+1$. If $r\ne \ell+1$, then the canonical point on level $\ell+1$ would occur in $\Theta$ but would have to lie both above the first $\ell-1$ ranks and below rank $\ell$, contradicting $\Theta^{-1}(\ell)_y>\ell+1$. Therefore $r=\ell+1$. Matching the remaining canonical levels in increasing order yields \eqref{eq:y-shift}, and the position of rank $\ell$ between levels $\ell+1$ and $\ell+2$ gives
\(\ell+1<\Theta^{-1}(\ell)_y<\ell+2\).
\end{proof}

We now derive the final contradiction. From \eqref{eq:x-shift},
\(\Theta(1)_y>\Theta(2)_y=\Sigma(1)_y=\sigma(1)\).
Because $\Theta$ is order-isomorphic to $\sigma$, this implies $\sigma(1)>\sigma(2)$. If $1\le i\le u-3$ and $\sigma(i)>\sigma(i+1)$, then \eqref{eq:x-shift} gives
\(\Theta(i+1)=\Sigma(i)\), \(\Theta(i+2)=\Sigma(i+1)\),
so the same order comparison in $\Theta$ implies $\sigma(i+1)>\sigma(i+2)$. Induction gives
\inlinetag{\sigma(1)>\sigma(2)>\cdots>\sigma(u-1)}{4.10}{eq:left-decreasing}.
Similarly, $\Theta(u)=\Theta^{-1}(k)$ lies above the entire canonical square, while $\Theta(u+1)=\Sigma(u)=\Sigma^{-1}(k)$. Hence $\sigma(u)>\sigma(u+1)$, and the same shift argument gives
\inlinetag{\sigma(u)>\sigma(u+1)>\cdots>\sigma(k)}{4.11}{eq:right-decreasing}.
Since $\sigma(k)=r=\ell+1>\ell=\sigma(1)$, every value in the second decreasing block is larger than every value in the first. Therefore
\(\sigma=(\ell(\ell-1)\cdots21)(k(k-1)\cdots(r+1)r)\), \(r=\ell+1\),
which is precisely one of the patterns treated in Subsection~\ref{subsec:two-blocks} and is excluded by the standing assumptions. This contradiction proves the lemma.
\end{proof}

\subsubsection{Both boundary gaps have size at least two}\label{subsubsec:large-gaps}

If $\ell<r$ and $r-\ell\ge2$, then $\PR$ meets every horizontal stripe $H_i$ with $i\le r-2$, while $\PL$ meets every $H_i$ with $i\ge \ell+1$; these two ranges cover $[k-1]$. The case $\ell>r$ is symmetric. Likewise, if $|b-u|\ge2$, the sets $\PB$ and $\PU$ together meet every vertical stripe.

\begin{corollary}\label{cor:large-gaps}
Theorem~\ref{thm:geometric-construction} holds whenever
\(|\ell-r|\ge2\) and \(|b-u|\ge2\).
\end{corollary}

\begin{proof}
Construction~\ref{con:four-boundary} covers every horizontal and vertical stripe under these two inequalities, and Lemma~\ref{lem:four-boundary-unique} gives uniqueness. The discretization principle preceding Subsection~\ref{subsec:two-blocks} then yields Theorem~\ref{thm:geometric-construction}.
\end{proof}

\subsubsection{At least one boundary gap has size one}\label{subsubsec:unit-gap}

It remains to handle the case $|\ell-r|=1$ or $|b-u|=1$. By a dihedral symmetry, we may assume
\inlinetag{\ell=r-1}{4.12}{eq:l-r-adjacent}.
Under this assumption, Construction~\ref{con:four-boundary} covers every horizontal stripe except $H_\ell$.

\begin{construction}[Extended boundary construction]\label{con:one-internal-segment}
For each integer $a$, let
\(P_a=(a,\ell)+t(1,1)\), \(0<t<1\).
Thus $P_a\subset H_\ell\cap V_a$. For a fixed $a\in\{2,3,\ldots,k-2\}$, write
\(D_a=\Sigma\cup\PL\cup\PR\cup\PB\cup\PU\cup P_a\).
\end{construction}

\begin{figure}[H]
\centering
\begin{tikzpicture}[scale=.48]
 \def\K{8}
 \draw[step=1,gray!42,thin] (1,1) grid (8,8);
 \draw[thick] (1,1) rectangle (8,8);
 \fill[red!10] (1,3) rectangle (8,4);
 \node[red!70!black,right] at (8.1,3.5) {$H_\ell$};
 \draw[blue!70!black,densely dotted,very thick] (1,4)--(-3,8);
 \draw[blue!70!black,densely dotted,very thick] (8,3)--(12,-1);
 \draw[blue!70!black,densely dotted,very thick] (3,1)--(8,-4);
 \draw[blue!70!black,densely dotted,very thick] (6,8)--(1,13);
 \draw[green!45!black,very thick] (5,3)--(6,4);
 \node[green!35!black,above] at (5.5,3.95) {$P_a$};
\end{tikzpicture}
\caption{Construction~\ref{con:one-internal-segment}: the added segment $P_a$ covers the single horizontal stripe omitted by the boundary construction.}
\end{figure}
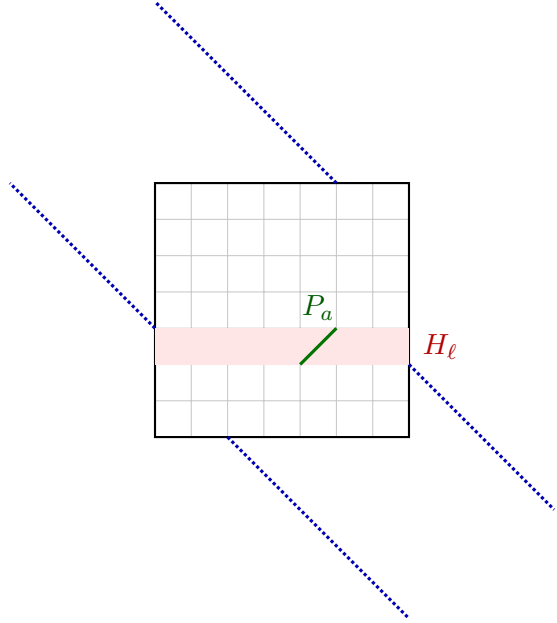

\begin{lemma}\label{lem:one-internal-segment}
There exists an index $a\in\{2,3,\ldots,k-2\}$ such that every point set $\Pi$ associated with a permutation and satisfying
\(\Sigma\subseteq\Pi\subseteq D_a\)
has $\Sigma$ as its unique $\sigma$-copy.
\end{lemma}

\begin{remark}
The restriction $2\le a\le k-2$ is essential for the argument. If $a\le1$ or $a\ge k-1$, the segment lies next to an outer boundary and can replace a boundary point of the canonical copy without forcing the propagation contradictions used below.
\end{remark}

\begin{proof}
Fix $a\in\{2,\ldots,k-2\}$ and suppose that $\Theta\subseteq D_a$ is a $\sigma$-copy distinct from $\Sigma$. By Lemma~\ref{lem:four-boundary-unique}, $\Theta$ must contain a point of $P_a$.

The boundary restrictions \eqref{eq:boundary-restrictions} continue to hold. We first adapt Claim~\ref{clm:propagation} to account for the internal segment.

\begin{claim}\label{clm:propagation-with-Pa}
The following statements hold.
\begin{enumerate}[label=\arabic*.,leftmargin=2.2em]
\item If $j<k$, $\Theta^{-1}(j)\notin P_a$, and $\Theta^{-1}(j)_y>j$, then $\Theta^{-1}(j+1)_y>j+1$.
\item If $j>1$, $\Theta^{-1}(j)\notin P_a$, and $\Theta^{-1}(j)_y<j$, then $\Theta^{-1}(j-1)_y<j-1$.
\item If $i<k$, $\Theta(i)\notin P_a$, and $\Theta(i)_x>i$, then $\Theta(i+1)_x>i+1$.
\item If $i>1$, $\Theta(i)\notin P_a$, and $\Theta(i)_x<i$, then $\Theta(i-1)_x<i-1$.
\end{enumerate}
\end{claim}

\begin{proof}
Once the current point is known not to lie in $P_a$, the proof of Claim~\ref{clm:propagation} applies verbatim.
\end{proof}

\begin{corollary}\label{cor:propagation-with-Pa}
The following statements hold.
\begin{enumerate}[label=\arabic*.,leftmargin=2.2em]
\item If $\Theta^{-1}(j)_y>j$ and $\Theta^{-1}(j)_y\ge \ell+1$, then
\(\Theta^{-1}(s)_y>s\), \((j\le s\le k)\).
\item If $\Theta^{-1}(j)_y<j$ and $\Theta^{-1}(j)_y\le \ell$, then
\(\Theta^{-1}(s)_y<s\), \((1\le s\le j)\).
\item If $\Theta(i)_x>i$ and $\Theta(i)_x\ge a+1$, then
\(\Theta(s)_x>s\), \((i\le s\le k)\).
\item If $\Theta(i)_x<i$ and $\Theta(i)_x\le a$, then
\(\Theta(s)_x<s\), \((1\le s\le i)\).
\end{enumerate}
\end{corollary}

\begin{proof}
We prove statement~1. The segment $P_a$ has $y$-coordinate range $(\ell,\ell+1)$. Hence a point with $y\ge \ell+1$ does not lie in $P_a$. Starting from $\Theta^{-1}(j)$, statement~1 of Claim~\ref{clm:propagation-with-Pa} applies; the next point has still larger $y$-coordinate and therefore also lies outside $P_a$. Iterating proves the conclusion. The proofs of statements~2--4 are identical, using respectively the ranges $y\le \ell$, $x\ge a+1$, and $x\le a$ to exclude $P_a$ at every step.
\end{proof}

\begin{claim}\label{clm:endpoints-not-Pa}
We have
\(\Theta(1)\notin P_a\) and \(\Theta(k)\notin P_a\).
\end{claim}

\begin{proof}
The only support points with $y$-coordinate in the open strip $(\ell,\ell+1)$ are the points of $P_a$. Consequently, the points of $\Theta\cap P_a$ occupy a consecutive interval in the value-rank order of $\Theta$. Since $P_a$ has positive slope, their position order agrees with their value-rank order.

Suppose $\Theta(1)\in P_a$. Its value-rank is $\ell$. Since $P_a\subset(a,a+1)\times(\ell,\ell+1)$ and $a\ge2$, we have $\Theta(1)_x>a\ge2$, and as $\Theta(1)$ is the leftmost point of $\Theta$, every other point of $\Theta$ has $x>2$; in particular $\Sigma(1)=(1,\ell)\notin\Theta$. Assume first that $\Theta(1)$ is the only point of $\Theta$ on $P_a$, and set $\Theta' = (\Theta\setminus \{\Theta(1)\}) \cup \{\Sigma(1)\}$. Then $|\Theta'|=k$, and $\Theta'$ lies in the support $\Sigma\cup\PL\cup\PR\cup\PB\cup\PU$ of Construction~\ref{con:four-boundary}. It is a $\sigma$-copy, as every selected point of value-rank below $\ell$ lies below the level $y=\ell$, and every selected point of value-rank above $\ell$ lies above that level, because there are no other support points in $(\ell,\ell+1)$; and $\Sigma(1)$ is strictly to the left of every other point of $\Theta'$. Finally $\Theta'\ne\Sigma$, because every point of $\Theta'$ other than $\Sigma(1)$ has $x>2$, so that $\Sigma(2)\in\Sigma\setminus\Theta'$. This contradicts Lemma~\ref{lem:four-boundary-unique}.

If $\Theta$ contains another point of $P_a$, then the next value-rank represented on $P_a$ must be $\ell+1=r$. The unique point of that rank is $\Theta(k)$; because $\Theta(k)$ is the last point in the position order, no further point of $P_a$ can be selected. Thus
\(\Theta\cap P_a=\{\Theta(1),\Theta(k)\}\).
Every intermediate point $\Theta(2),\ldots,\Theta(k-1)$ has $x$-coordinate strictly between the $x$-coordinates of these two points, hence in $(a,a+1)$. It cannot be canonical, and it cannot lie in $\PL$ or $\PR$; therefore it lies in $\PB\cup\PU$. The boundary restrictions allow at most one selected point in each of $\PB$ and $\PU$, so $k-2\le2$. The four boundary points are distinct, so $k\ge4$, and hence $k=4$. Then $a=2$, $\ell=2$, and $r=3$. The only two possible patterns are $(2143)$ and $(2413)$. For $(2143)$, one has $b<u$ and Construction~\ref{con:four-boundary} places $\PU$ in $x<2$ and $\PB$ in $x>3$; for $(2413)$, one has $b>u$ and the two inequalities are reversed. In neither case can an intermediate point lie in $(2,3)$, a contradiction. Hence $\Theta(1)\notin P_a$.

The argument for $\Theta(k)$ is the same, read from the other end. If it is the only selected point on $P_a$, replacing it by $\Sigma(k)=(k,r)$ gives a copy in the boundary-only support. If another point of $P_a$ is selected, consecutiveness of the represented value-ranks forces that point to be $\Theta(1)$, and the preceding $k=4$ argument applies. Therefore $\Theta(k)\notin P_a$ as well.
\end{proof}

Since $\Theta$ contains a point of $P_a$ and $\Theta^{-1}(\ell)=\Theta(1)$ and $\Theta^{-1}(r)=\Theta(k)$ do not lie there, some value-rank $s$ with either $s<\ell$ or $s>r$ satisfies $\Theta^{-1}(s)\in P_a$.

\begin{claim}\label{clm:shift-from-Pa}
The following statements hold.
\begin{enumerate}[label=\arabic*.,leftmargin=2.2em]
\item If $\Theta^{-1}(s)\in P_a$ for some $s<\ell$, then
\(\Theta^{-1}(\ell)_y\ge \ell+1\), \(\Theta^{-1}(k)\in\PU\),
and
\inlinetag{\Theta^{-1}(j)=\Sigma^{-1}(j+1) \quad\allowbreak(r\le j\le k-1)}{4.13}{eq:upper-rank-shift}.
\item If $\Theta^{-1}(s)\in P_a$ for some $s>r$, then
\(\Theta^{-1}(r)_y\le \ell\), \(\Theta^{-1}(1)\in\PB\),
and
\inlinetag{\Theta^{-1}(j)=\Sigma^{-1}(j-1) \quad\allowbreak(2\le j\le \ell)}{4.14}{eq:lower-rank-shift}.
\end{enumerate}
\end{claim}

\begin{proof}
We prove statement~1. Since $P_a$ has $y$-coordinate range $(\ell,\ell+1)$ and $s<\ell$,
\(\Theta^{-1}(\ell)_y>\Theta^{-1}(s)_y>\ell\).
Claim~\ref{clm:endpoints-not-Pa} gives $\Theta^{-1}(\ell)=\Theta(1)\notin P_a$. No other part of the support has a point with $y\in(\ell,\ell+1)$, so $\Theta^{-1}(\ell)_y\ge \ell+1$. Corollary~\ref{cor:propagation-with-Pa} now gives $\Theta^{-1}(k)_y>k$, hence $\Theta^{-1}(k)\in\PU$.

For $r\le j\le k-1$, the point $\Theta^{-1}(j)$ has $y>r$. It cannot lie in $P_a$ or $\PR$, whose $y$-coordinates are below $r$, and it is excluded from $\PL$, $\PB$, and $\PU$ by the boundary restrictions because its rank and position are not the corresponding boundary labels. Thus these $k-r$ points lie in $\Sigma$. The only $k-r$ canonical points with $y>r$ are
\(\Sigma^{-1}(r+1),\ldots,\Sigma^{-1}(k)\).
Matching in increasing $y$-order proves \eqref{eq:upper-rank-shift}. Statement~2 follows similarly in the reverse vertical order.
\end{proof}

Define the two monotonicity conditions
\inlinetag{\sigma^{-1}(k)<\sigma^{-1}(k-1)<\cdots<\sigma^{-1}(r+1)<\sigma^{-1}(r)}{UD}{eq:UD}
and
\inlinetag{\sigma^{-1}(\ell)<\sigma^{-1}(\ell-1)<\cdots<\sigma^{-1}(2)<\sigma^{-1}(1)}{BD}{eq:BD}.
\begin{claim}\label{clm:monotonicity-forced}
The following statements hold.
\begin{enumerate}[label=\arabic*.,leftmargin=2.2em]
\item If $\Theta^{-1}(s)\in P_a$ for some $s<\ell$, then \eqref{eq:UD} holds,
\(\Theta^{-1}(k)_x<u-1\),
and $b<u$.
\item If $\Theta^{-1}(s)\in P_a$ for some $s>r$, then \eqref{eq:BD} holds,
\(\Theta^{-1}(1)_x>b+1\),
and $b<u$.
\end{enumerate}
\end{claim}

\begin{proof}
We prove statement~1. By \eqref{eq:upper-rank-shift},
\(\Theta^{-1}(j)=\Sigma^{-1}(j+1)\), \((r\le j\le k-1)\).
Since $\sigma^{-1}(r)=k$, the relative $x$-order of the points of ranks $r,r+1$ in $\Theta$ forces
\(\sigma^{-1}(r+2)<\sigma^{-1}(r+1)<\sigma^{-1}(r)\) (if \(r\le k-2)\), or \(\sigma^{-1}(r+1)<\sigma^{-1}(r)\) (if \(r = k-1)\).
Repeating the same comparison for consecutive ranks yields \eqref{eq:UD}.

In particular,
\(\Theta^{-1}(k)_x<\Theta^{-1}(k-1)_x =\Sigma^{-1}(k)_x=u\).
If $b>u$, every point of $\PU$ has $x>u+1$, contradicting $\Theta^{-1}(k)\in\PU$. Hence $b<u$, and then the trimming of $\PU$ gives $\Theta^{-1}(k)_x<u-1$. Statement~2 follows similarly in the reverse vertical order.
\end{proof}

\begin{corollary}\label{cor:generic-a}
If $b>u$, or if neither \eqref{eq:UD} nor \eqref{eq:BD} holds, then every $a\in\{2,\ldots,k-2\}$ satisfies the conclusion of Lemma~\ref{lem:one-internal-segment}.
\end{corollary}

\begin{proof}
Any additional copy must contain a point of $P_a$ whose value-rank lies below $\ell$ or above $r$. Claim~\ref{clm:monotonicity-forced} then forces $b<u$ and, respectively, \eqref{eq:UD} or \eqref{eq:BD}. Under either hypothesis of the corollary this is impossible.
\end{proof}

If both \eqref{eq:UD} and \eqref{eq:BD} hold, then, together with $b<u$ and $r=\ell+1$, the positions of all values are forced into the order
\(\sigma^{-1}(\ell)<\cdots<\sigma^{-1}(1)< \sigma^{-1}(k)<\cdots<\sigma^{-1}(r)\),
so $\sigma$ is the two-block pattern treated in Subsection~\ref{subsec:two-blocks}. We may therefore assume that exactly one of \eqref{eq:UD} and \eqref{eq:BD} holds. By rotational symmetry, it is enough to treat the case in which \eqref{eq:UD} holds and \eqref{eq:BD} fails. Then every additional copy must have $\Theta^{-1}(s)\in P_a$ for some $s<\ell$.

\begin{claim}\label{clm:contiguous-Pa-ranks}
If $\Theta^{-1}(s)\in P_a$ for some $s<\ell$, then
\(\Theta^{-1}(t)\in P_a\), \((s\le t<\ell)\).
\end{claim}

\begin{proof}
Fix $t$ with $s\le t<\ell$. Since $\Theta^{-1}(s)\in P_a$,
\(\ell<\Theta^{-1}(t)_y\).
By \eqref{eq:upper-rank-shift},
\(\Theta^{-1}(r)=\Sigma^{-1}(r+1)\),
so $\Theta^{-1}(r)_y=r+1=\ell+2$. As $t<r$,
\(\ell<\Theta^{-1}(t)_y<\ell+2\).
The only parts of $D_a$ meeting this open horizontal band are $P_a$, $\PL$, and the canonical point $\Sigma^{-1}(r)$. The point $\Theta^{-1}(t)$ is not in $\PL$, because $t<\ell$ while the boundary restriction allows only $\Theta^{-1}(\ell)=\Theta(1)$ there. It is not $\Sigma^{-1}(r)$ either: since $\sigma^{-1}(r)=k$ and \eqref{eq:UD} gives $\sigma^{-1}(r+1)<\sigma^{-1}(r)$,
\(\Theta^{-1}(t)_x<\Theta^{-1}(r)_x =\Sigma^{-1}(r+1)_x <\Sigma^{-1}(r)_x\).
Therefore $\Theta^{-1}(t)\in P_a$.
\end{proof}

Return to the proof of Lemma~\ref{lem:one-internal-segment}. Claim~\ref{clm:contiguous-Pa-ranks} gives
\inlinetag{\Theta^{-1}(\ell-1)\in P_a, \quad\allowbreak a<\Theta^{-1}(\ell-1)_x<a+1}{4.15}{eq:lminus1-in-Pa}.
Claim~\ref{clm:shift-from-Pa} also gives
\inlinetag{\Theta^{-1}(k-1)=\Sigma^{-1}(k), \quad\allowbreak \Theta^{-1}(k-1)_x=u}{4.16}{eq:kminus1-at-u},
while Claim~\ref{clm:monotonicity-forced} gives
\inlinetag{\Theta^{-1}(k)_x<u-1}{4.17}{eq:top-left-u}.
Let
\(p=\sigma^{-1}(\ell-1)\), \(q=\sigma^{-1}(k-1)\).
If $p<u$ or $p>q$, choose $a=u-1$. This is a valid index: $b<u$ and $b\ge2$ imply $u\ge3$, while \eqref{eq:UD} gives $u<q\le k$, hence $u\le k-1$. Equations \eqref{eq:lminus1-in-Pa}--\eqref{eq:top-left-u} then yield
\(\Theta^{-1}(k)_x<u-1<\Theta^{-1}(\ell-1)_x<u =\Theta^{-1}(k-1)_x\).
Because $\Theta$ is order-isomorphic to $\sigma$, this forces
\(u<p<q\),
contradicting the present case.

It remains that
\(u<p<q\).
Choose $a=k-2$. Since $u<p<q\le k$, we have $u\le k-2$. Equations \eqref{eq:lminus1-in-Pa} and \eqref{eq:kminus1-at-u} give
\(\Theta^{-1}(k-1)_x=u\le k-2 <\Theta^{-1}(\ell-1)_x\).
Thus $q<p$, again a contradiction. In all cases we have produced a valid index $a$ for which no additional copy exists, proving the lemma.
\end{proof}

The preceding proof also identifies when the segment needed for uniqueness can be chosen to cover a second missing stripe.

\begin{corollary}\label{cor:horizontal-exceptions}
Assume $\ell=r-1$. Then Theorem~\ref{thm:geometric-construction} holds unless $b=u-1$ and one of the following two configurations occurs:
\begin{enumerate}[label=(\alph*),leftmargin=2.2em]
\item \eqref{eq:UD} holds, \eqref{eq:BD} fails, and
\inlinetag{u<\sigma^{-1}(\ell-1)<\sigma^{-1}(k-1)}{4.18}{eq:horizontal-exception-a};
\item \eqref{eq:BD} holds, \eqref{eq:UD} fails, and
\inlinetag{\sigma^{-1}(2)<\sigma^{-1}(r+1)<b}{4.19}{eq:horizontal-exception-b}.
\end{enumerate}
\end{corollary}

\begin{proof}
If $|b-u|\ge2$, then the boundary sets $\PB$ and $\PU$ already cover every vertical stripe, and Lemma~\ref{lem:one-internal-segment} supplies a segment $P_a$ that covers the missing horizontal stripe $H_\ell$ while preserving uniqueness. If $b=u+1$, then the only missing vertical stripe is $V_u$; because $b>u$, Corollary~\ref{cor:generic-a} shows that every admissible segment preserves uniqueness, so we take $a=u$ and cover $H_\ell$ and $V_u$ simultaneously.

It remains to assume $b=u-1$. Then the boundary construction misses exactly the vertical stripe $V_b$, so the added segment must be $P_b=P_{u-1}$ in order to cover both missing stripes. The proof of Lemma~\ref{lem:one-internal-segment} shows that this choice works whenever neither exceptional ordering in \eqref{eq:horizontal-exception-a}--\eqref{eq:horizontal-exception-b} occurs. If both \eqref{eq:UD} and \eqref{eq:BD} hold, the pattern belongs to Subsection~\ref{subsec:two-blocks} and has already been handled. This proves the corollary.
\end{proof}

For the coordinate-dual statement, define
\inlinetag{\sigma(1)>\sigma(2)>\cdots>\sigma(b-1)>\sigma(b)}{LD}{eq:LD}
and
\inlinetag{\sigma(u)>\sigma(u+1)>\cdots>\sigma(k-1)>\sigma(k)}{RD}{eq:RD}.
\begin{corollary}\label{cor:vertical-exceptions}
Assume $b=u-1$. Then Theorem~\ref{thm:geometric-construction} holds unless $\ell=r-1$ and one of the following two configurations occurs:
\begin{enumerate}[label=(\alph*),leftmargin=2.2em]
\item \eqref{eq:LD} holds, \eqref{eq:RD} fails, and
\inlinetag{\sigma(2)<\sigma(u+1)<\ell}{4.20}{eq:vertical-exception-a};
\item \eqref{eq:RD} holds, \eqref{eq:LD} fails, and
\inlinetag{r<\sigma(b-1)<\sigma(k-1)}{4.21}{eq:vertical-exception-b}.
\end{enumerate}
\end{corollary}

\begin{proof}
Apply Corollary~\ref{cor:horizontal-exceptions} after reflecting the configuration in
the anti-diagonal via
\(A(x,y)=(k+1-y,k+1-x)\).
Under this reflection $\ell$ and $u$ exchange roles, as do $r$ and $b$; the hypothesis
$b=u-1$ becomes $\ell=r-1$; the conditions \eqref{eq:UD} and \eqref{eq:BD} become
\eqref{eq:LD} and \eqref{eq:RD}; and the two exceptional inequalities become
\eqref{eq:vertical-exception-a} and \eqref{eq:vertical-exception-b}.
\end{proof}

\subsubsection{The simultaneous exceptional case}\label{subsubsec:simultaneous-exception}

It remains to consider patterns satisfying an exceptional configuration in both Corollaries~\ref{cor:horizontal-exceptions} and~\ref{cor:vertical-exceptions}. Without loss of generality, we may assume that \eqref{eq:horizontal-exception-a} holds. Thus
\inlinetag{\ell=r-1, \quad\allowbreak b=u-1}{4.22}{eq:adjacent-both},
\eqref{eq:UD} holds, \eqref{eq:BD} fails, and
\inlinetag{u<\sigma^{-1}(\ell-1)<\sigma^{-1}(k-1)}{4.23}{eq:final-position-condition}.
Under these assumptions, \eqref{eq:RD} cannot hold. Indeed, the two positions
\(\sigma^{-1}(\ell-1)<\sigma^{-1}(k-1)\)
both lie to the right of $u$, but their values satisfy $\ell-1<k-1$, contradicting the decrease required by \eqref{eq:RD}. Consequently the exceptional case in Corollary~\ref{cor:vertical-exceptions} must be \eqref{eq:vertical-exception-a}. Hence \eqref{eq:LD} holds and
\inlinetag{\sigma(2)<\sigma(u+1)<\ell}{4.24}{eq:final-value-condition}.
\begin{construction}[Two internal segments]\label{con:two-internal-segments}
Starting from Construction~\ref{con:four-boundary}, add
\(\PH=(k-2,\ell)+t(1,1)\), \(0<t<1\),
and
\(\PV=(u-1,2)+t(1,1)\), \(0<t<1\).
The segment $\PH$ lies in the missing horizontal stripe $H_\ell$, and $\PV$ lies in the missing vertical stripe $V_{u-1}=V_b$.
\end{construction}

\begin{figure}[H]
\centering
\begin{tikzpicture}[scale=.47]
 \def\K{8}
 \draw[step=1,gray!42,thin] (1,1) grid (8,8);
 \draw[thick] (1,1) rectangle (8,8);
 \fill[red!10] (1,5) rectangle (8,6);
 \fill[red!10] (3,1) rectangle (4,8);
 \fill[red!75!black] (1,5) circle (2.4pt) node[left=2pt] {$L$};
\fill[red!75!black] (8,6) circle (2.4pt) node[right=2pt] {$R$};
\fill[red!75!black] (3,1) circle (2.4pt) node[below=2pt] {$B$};
\fill[red!75!black] (4,8) circle (2.4pt) node[above=2pt] {$U$};
 \draw[blue!70!black,densely dotted,very thick] (1,8)--(-1,6);
 \draw[blue!70!black,densely dotted,very thick] (8,5)--(12,1);
 \draw[blue!70!black,densely dotted,very thick] (4,1)--(8,-3);
 \draw[blue!70!black,densely dotted,very thick] (1,8)--(3,10);
 \draw[green!45!black,very thick] (6,5)--(7,6);
 \draw[green!45!black,very thick] (3,2)--(4,3);
 \node[green!35!black,above] at (6.5,5.95) {$P_{\mathrm H}$};
 \node[green!35!black,left] at (3.35,2.45) {$P_{\mathrm V}$};
 \draw[blue!70!black,fill=white,thin] (1,8) circle (2.2pt);
\end{tikzpicture}
\caption{Construction~\ref{con:two-internal-segments}. The two internal segments cover the horizontal and vertical stripes left uncovered by the boundary sets.}
\end{figure}
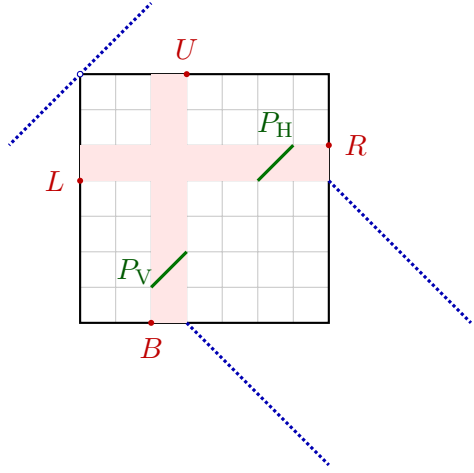

\begin{lemma}\label{lem:two-internal-segments}
Assume the standing hypotheses of Subsection~\ref{subsec:four-boundary-points}, together with \eqref{eq:adjacent-both}--\eqref{eq:final-value-condition}, \eqref{eq:UD}, and \eqref{eq:LD}. If $\Pi$ is associated with a permutation and
\(\Sigma\subseteq\Pi\subseteq \Sigma\cup\PL\cup\PR\cup\PB\cup\PU\cup\PH\cup\PV\),
then $\Sigma$ is the unique $\sigma$-copy in $\Pi$.
\end{lemma}

\begin{remark}
The numerical conditions of Lemma~\ref{lem:two-internal-segments} are satisfied by the pattern $(31524)$, and for that pattern the support of Construction~\ref{con:two-internal-segments} may contain additional copies. This does not affect the lemma, as the standing hypotheses of Subsection~\ref{subsec:four-boundary-points} exclude every pattern already treated in Subsection~\ref{subsec:exceptional-five}, so $(31524)$ does not satisfy the hypotheses of Lemma~\ref{lem:two-internal-segments}. The proof below shows that $(31524)$ is the only excluded pattern that can arise from the final reduction.
\end{remark}

\begin{proof}
Suppose, for contradiction, that $\Theta$ is a $\sigma$-copy in the support of Construction~\ref{con:two-internal-segments} with $\Theta\ne\Sigma$. The support obtained by deleting $\PV$ is the choice $a=k-2$ in Construction~\ref{con:one-internal-segment}; the proof of Lemma~\ref{lem:one-internal-segment}, in the exceptional ordering \eqref{eq:final-position-condition}, shows that it has no additional copy. The coordinate-dual statement applies after deleting $\PH$. Therefore
\inlinetag{\Theta\cap\PH\ne\varnothing \quad\allowbreak\text{and}\quad\allowbreak \Theta\cap\PV\ne\varnothing}{4.25}{eq:uses-both-internal}.
The inequalities in \eqref{eq:final-position-condition} imply $u\le k-2$, while \eqref{eq:final-value-condition} implies $\ell\ge3$. Hence
\inlinetag{\PH\subset(k-2,k-1)\times(\ell,\ell+1), \quad\allowbreak \PV\subset(u-1,u)\times(2,3)}{4.26}{eq:internal-ranges},
and these two segments have disjoint $x$- and $y$-coordinate ranges.

The boundary restrictions \eqref{eq:boundary-restrictions} remain valid. We next record the propagation statements needed in the enlarged support.

\begin{claim}\label{clm:propagation-with-two}
The following statements hold.
\begin{enumerate}[label=\arabic*.,leftmargin=2.2em]
\item If $j<k$, $\Theta^{-1}(j)\notin\PH\cup\PV$, and $\Theta^{-1}(j)_y>j$, then $\Theta^{-1}(j+1)_y>j+1$.
\item If $j>1$, $\Theta^{-1}(j)\notin\PH\cup\PV$, and $\Theta^{-1}(j)_y<j$, then $\Theta^{-1}(j-1)_y<j-1$.
\item If $i<k$, $\Theta(i)\notin\PH\cup\PV$, and $\Theta(i)_x>i$, then $\Theta(i+1)_x>i+1$.
\item If $i>1$, $\Theta(i)\notin\PH\cup\PV$, and $\Theta(i)_x<i$, then $\Theta(i-1)_x<i-1$.
\end{enumerate}
\end{claim}

\begin{proof}
The proof is the same as for Claim~\ref{clm:propagation-with-Pa}: after excluding the two internal segments, Claim~\ref{clm:propagation} applies.
\end{proof}

The following corollary is the analogue of Corollaries~\ref{cor:propagation}
and~\ref{cor:propagation-with-Pa} under Construction~\ref{con:two-internal-segments}.
Unlike the earlier setting, the present construction is not symmetric with respect to
all four directions. Due to the coordinate ranges of $\PH$ and $\PV$, only the upward
propagation along the positive $y$-direction and the leftward propagation in the
negative $x$-direction remain valid in the form stated below.

\begin{corollary}\label{cor:one-sided-propagation-two}
The following statements hold.
\begin{enumerate}[label=\arabic*.,leftmargin=2.2em]
\item If $\Theta^{-1}(j)_y>j$ and $\Theta^{-1}(j)_y\ge \ell+1$, then
\(\Theta^{-1}(s)_y>s\), \((j\le s\le k)\).
\item If $\Theta(i)_x<i$ and $\Theta(i)_x\le u-1$, then
\(\Theta(s)_x<s\), \((1\le s\le i)\).
\end{enumerate}
\end{corollary}

\begin{proof}
For statement~1, a point with $y\ge \ell+1$ lies outside both internal segments: $\PH$ has $y<\ell+1$, and $\PV$ has $y<3\le \ell$. Thus Claim~\ref{clm:propagation-with-two} can be applied repeatedly. For statement~2, a point with $x\le u-1$ lies outside $\PV$, whose $x$-coordinates exceed $u-1$, and outside $\PH$, whose $x$-coordinates exceed $k-2\ge u$. Repeatedly apply statement~4 of Claim~\ref{clm:propagation-with-two}.
\end{proof}

\begin{claim}\label{clm:boundary-ranks-not-internal}
We have
\(\Theta^{-1}(\ell)\notin\PH\) and \(\Theta(u)\notin\PV\).
\end{claim}

\begin{proof}
Suppose first that $\Theta^{-1}(\ell)=\Theta(1)\in\PH$. Then
\(k-2<\Theta^{-1}(\ell)_x<k-1\), \(\ell<\Theta^{-1}(\ell)_y<\ell+1\).
The top-rank point $\Theta^{-1}(k)=\Theta(u)$ lies strictly northeast of it. By \eqref{eq:internal-ranges} and the trimming of the boundary sets, the only support points with both $x>k-2$ and $y>\ell$ lie in
\(\PH\cup\{\Sigma(k-1),\Sigma(k)\}\).

If $\Theta^{-1}(k)\in\PH$, then the intermediate rank-$r$ point $\Theta^{-1}(r)=\Theta(k)$ also lies in $\PH$. In the $x$-order the three relevant ranks occur as
\(\Theta^{-1}(\ell)=\Theta(1),\quad \Theta^{-1}(k)=\Theta(u),\quad \Theta^{-1}(r)=\Theta(k)\),
because $1<u<k$, whereas their $y$-ranks are ordered $\ell<r<k$. Thus the last two points form an inversion on the positive-slope segment $\PH$, a contradiction.
If $\Theta^{-1}(k)=\Sigma(k)$, then $\Theta(k)$ has $y$-coordinate strictly between $\ell$ and $\ell+1$ and hence lies in $\PH$, but its $x$-coordinate must exceed the $x$-coordinate $k$ of $\Theta^{-1}(k)$, an impossibility.

The remaining possibility is $\Theta^{-1}(k)=\Sigma(k-1)$. Then $\Theta^{-1}(k)_x=k-1$ and $\Theta^{-1}(k)_y=\sigma(k-1)>\ell$. Since $u\le k-2$, both $\Theta(k-1)$ and $\Theta(k)$ lie to the right of $\Theta(u)=\Theta^{-1}(k)$; their values $\sigma(k-1)$ and $r$ are both larger than $\ell$. Thus both points would have $x>k-1$ and $y>\ell$. The support contains at most the single point $\Sigma(k)$ with these two inequalities, a contradiction. Hence $\Theta^{-1}(\ell)\notin\PH$.

For the second assertion, reflect the entire configuration in the anti-diagonal via
\(A(x,y)=(k+1-y,k+1-x)\).
The transformed pattern has boundary parameters
\(\ell'=k+1-u\), \(r'=\ell'+1\), \(b'=k+1-r\), \(u'=b'+1\).
Moreover, $A(\PV)$ is the horizontal internal segment
\((k-2,\ell')+t(1,1)\), \(0<t<1\),
while $A(\PH)$ is the corresponding vertical internal segment. Under the same reflection, \eqref{eq:LD} and \eqref{eq:final-value-condition} become the analogues of \eqref{eq:UD} and \eqref{eq:final-position-condition}. Applying the first assertion to the reflected copy and translating back gives $\Theta(u)\notin\PV$.
\end{proof}

We can now repeat the rank-shift arguments from Claims~\ref{clm:shift-from-Pa} and~\ref{clm:contiguous-Pa-ranks}; the only places where the earlier proof used the absence of another internal segment are supplied here by Corollary~\ref{cor:one-sided-propagation-two} and Claim~\ref{clm:boundary-ranks-not-internal}. For completeness, we state the resulting conclusions and indicate the key exclusions.

\begin{claim}\label{clm:shift-from-PH-PV}
The following statements hold.
\begin{enumerate}[label=\arabic*.,leftmargin=2.2em]
\item If $\Theta^{-1}(s)\in\PH$ for some $s<\ell$, then
\(\Theta^{-1}(\ell)_y\ge \ell+1\), \(\Theta^{-1}(k)\in\PU\),
and
\inlinetag{\Theta^{-1}(j)=\Sigma^{-1}(j+1) \quad\allowbreak(r\le j\le k-1)}{4.27}{eq:PH-shift}.
\item If $\Theta(s)\in\PV$ for some $s>u$, then
\(\Theta(u)_x\le u-1\), \(\Theta(1)\in\PL\),
and
\inlinetag{\Theta(i)=\Sigma(i-1) \quad\allowbreak(2\le i\le b)}{4.28}{eq:PV-shift}.
\end{enumerate}
\end{claim}

\begin{proof}
For statement~1, Claim~\ref{clm:boundary-ranks-not-internal} replaces Claim~\ref{clm:endpoints-not-Pa} in the proof of Claim~\ref{clm:shift-from-Pa}. Once $\Theta^{-1}(\ell)_y\ge \ell+1$, Corollary~\ref{cor:one-sided-propagation-two} propagates upward. The points of ranks $r,\ldots,k-1$ lie above both internal segments, outside the trimmed right boundary, and are not boundary-labeled points; hence they lie in $\Sigma$ and match as in \eqref{eq:PH-shift}. Statement~2 is proved by the coordinate-dual argument, using the second part of Claim~\ref{clm:boundary-ranks-not-internal} and the leftward propagation in Corollary~\ref{cor:one-sided-propagation-two}.
\end{proof}

\begin{claim}\label{clm:contiguous-internal}
The following statements hold.
\begin{enumerate}[label=\arabic*.,leftmargin=2.2em]
\item If $\Theta^{-1}(s)\in\PH$ for some $s<\ell$, then
\(\Theta^{-1}(t)\in\PH\), \((s\le t<\ell)\).
\item If $\Theta(s)\in\PV$ for some $s>u$, then
\(\Theta(t)\in\PV\), \((u<t\le s)\).
\end{enumerate}
\end{claim}

\begin{proof}
For statement~1, the proof of Claim~\ref{clm:contiguous-Pa-ranks} applies with $\PH$ in place of $P_a$. The extra segment $\PV$ has $y$-coordinate range $(2,3)$, which is disjoint from the band $(\ell,\ell+2)$ because $\ell\ge3$. Statement~2 follows by the corresponding coordinate argument; the segment $\PH$ has $x>k-2\ge u$ and therefore does not interfere with the band used there.
\end{proof}

\begin{claim}\label{clm:which-ranks-use-internal}
We have
\(\Theta^{-1}(s)\notin\PH\quad(1\le s\le \ell)\)
and
\(\Theta(s)\notin\PV\quad(u\le s\le k)\).
\end{claim}

\begin{proof}
Claim~\ref{clm:boundary-ranks-not-internal} gives the first statement for $s=\ell$. If $\Theta^{-1}(s)\in\PH$ for some $s<\ell$, then Claim~\ref{clm:contiguous-internal} gives $\Theta^{-1}(\ell-1)\in\PH$, so
\(\Theta^{-1}(\ell-1)_x>k-2\).
Claim~\ref{clm:shift-from-PH-PV} gives
\(\Theta^{-1}(k-1)=\Sigma^{-1}(k)\), \(\Theta^{-1}(k-1)_x=u\le k-2\).
Thus the relative $x$-order in $\Theta$ implies
\(\sigma^{-1}(k-1)<\sigma^{-1}(\ell-1)\),
contradicting \eqref{eq:final-position-condition}. This proves the first assertion.

For the second assertion, Claim~\ref{clm:boundary-ranks-not-internal} covers $s=u$. If $\Theta(s)\in\PV$ for some $s>u$, then Claim~\ref{clm:contiguous-internal} gives $\Theta(u+1)\in\PV$, whose $y$-coordinate is less than $3$. Claim~\ref{clm:shift-from-PH-PV} gives $\Theta(2)=\Sigma(1)$, so $\Theta(2)_y=\ell\ge3$. This contradicts \eqref{eq:final-value-condition}, because $2<u+1$ and $\sigma(2)<\sigma(u+1)$ require $\Theta(2)_y<\Theta(u+1)_y$.
\end{proof}

By \eqref{eq:uses-both-internal} and Claim~\ref{clm:which-ranks-use-internal}, choose indices
\(h>\ell\) and \(v<u\)
such that
\inlinetag{\Theta^{-1}(h)\in\PH, \quad\allowbreak \Theta(v)\in\PV}{4.29}{eq:h-v}.
In particular,
\inlinetag{\ell<\Theta^{-1}(h)_y<\ell+1, \quad\allowbreak u-1<\Theta(v)_x<u}{4.30}{eq:h-v-ranges}.
Let
\(\Sigma_{\mathrm{BR}} = \Sigma\cap([u+1,k-1]\times[2,\ell-1])\)
and
\(\Theta_{\mathrm{BR}} = \{\Theta(i):\Sigma(i)\in\Sigma_{\mathrm{BR}}\}\).

\begin{claim}\label{clm:BR-contained}
We have
\(\Theta_{\mathrm{BR}}\subseteq\Sigma_{\mathrm{BR}}\).
\end{claim}

\begin{proof}
Fix $i$ with $\Sigma(i)\in\Sigma_{\mathrm{BR}}$. Thus
\(u+1\le i\le k-1\), \(2\le\sigma(i)\le \ell-1\).
Claim~\ref{clm:which-ranks-use-internal} excludes $\Theta(i)$ from $\PV$ and, because $\Theta(i)=\Theta^{-1}(\sigma(i))$, from $\PH$. Since neither its position nor its value-rank is a boundary label, \eqref{eq:boundary-restrictions} excludes all four boundary sets. Hence $\Theta(i)\in\Sigma$.

From \eqref{eq:h-v}, $v<i$ and $\sigma(i)<h$, so
\(\Theta(i)_x>\Theta(v)_x>u-1\), \(\Theta(i)_y<\Theta^{-1}(h)_y<\ell+1\).
Because $\Theta(i)\in\Sigma$, both coordinates are integers. Thus
\(u\le\Theta(i)_x\le k\), \(1\le\Theta(i)_y\le \ell\).
The values $\Theta(i)_x\in\{u,k\}$ are impossible, because $\Sigma(u)$ and $\Sigma(k)$ have $y$-coordinates $k$ and $r=\ell+1$, respectively. The values $\Theta(i)_y\in\{1,\ell\}$ are impossible, because $\Sigma^{-1}(1)$ and $\Sigma^{-1}(\ell)$ have $x$-coordinates $b=u-1$ and $1$, respectively. Therefore
\(u+1\le\Theta(i)_x\le k-1\), \(2\le\Theta(i)_y\le \ell-1\),
so $\Theta(i)\in\Sigma_{\mathrm{BR}}$.
\end{proof}

\begin{corollary}\label{cor:BR-fixed}
If
\(u+1\le i\le k-1\) and \(2\le\sigma(i)\le \ell-1\),
then
\(\Theta(i)=\Sigma(i)\).
\end{corollary}

\begin{proof}
By Claim~\ref{clm:BR-contained}, $\Theta_{\mathrm{BR}}$ is a subset of $\Sigma_{\mathrm{BR}}$. The two sets have the same cardinality by definition, so they are equal. Both are listed in increasing $x$-order by the same index set; hence corresponding points agree.
\end{proof}

\begin{claim}\label{clm:r-b-not-internal}
We have
\(\Theta^{-1}(r)\notin\PH\) and \(\Theta(b)\notin\PV\).
\end{claim}

\begin{proof}
Suppose $\Theta^{-1}(r)=\Theta(k)\in\PH$. Then
\(\Theta(k)_x<k-1\), \(\ell<\Theta(k)_y<\ell+1\).
The value $\sigma(k-1)$ is either at most $\ell-1$ or at least $r+1$.

If $\sigma(k-1)\le \ell-1$, then also $\sigma(k-1)\ge2$, otherwise $\sigma(k-1)=1$ would give $b=\sigma^{-1}(1)=k-1$, contradicting $b=u-1\le k-3$. Since moreover $u+1\le k-1$, the index $i=k-1$ satisfies the hypotheses of Corollary~\ref{cor:BR-fixed}, hence $\Theta(k-1)=\Sigma(k-1)$ and $\Theta(k-1)_x=k-1$, contradicting $\Theta(k-1)_x<\Theta(k)_x<k-1$.

Assume $\sigma(k-1)\ge r+1$. Condition \eqref{eq:UD} forces $\sigma(k-1)=r+1$: if it were larger, the value $r+1$ would have to occur strictly between positions $k-1$ and $k=\sigma^{-1}(r)$. Thus
\(\Theta(k-1)=\Theta^{-1}(r+1)\). Here $r+1<k$, otherwise $\sigma(k-1)=r+1=k$ would give $u=\sigma^{-1}(k)=k-1$, contradicting $u\le k-2$. Hence $\Theta(k-1)$ is none of the four points labeled in \eqref{eq:boundary-restrictions}, so it lies in none of the boundary sets. It does not lie in $\PV$, since its $y$-coordinate exceeds $\Theta(k)_y>\ell\ge3$, and it cannot lie in $\PH$, because $k-1<k$ but $r+1>r$ while $\PH$ has positive slope. Hence $\Theta(k-1)\in\Sigma$.

Moreover, $\Theta(k-1)_x<\Theta(k)_x<k-1$. The canonical levels $r$ and $r+1$ occur at $x=k$ and $x=k-1$, so $\Theta(k-1)_y>r+1$. Therefore
\(\Theta^{-1}(r+1)_y>r+1\).
Corollary~\ref{cor:one-sided-propagation-two} gives $\Theta^{-1}(k)_y>k$, and hence $\Theta^{-1}(k)\in\PU$. Since $b<u$, every point of $\PU$ has $x<u-1$. But \eqref{eq:h-v} gives $v<u$ and $\Theta(v)_x>u-1$, contradicting the required order $\Theta(v)_x<\Theta(u)_x=\Theta^{-1}(k)_x$. Thus $\Theta^{-1}(r)\notin\PH$.

For the assertion $\Theta(b)\notin\PV$, apply the same proof after the anti-diagonal reflection $A(x,y)=(k+1-y,k+1-x)$ used in Claim~\ref{clm:boundary-ranks-not-internal}. It interchanges the two internal segments, sends the bottom-right block to the corresponding reflected block, and exchanges the pairs \eqref{eq:UD}, \eqref{eq:final-position-condition} and \eqref{eq:LD}, \eqref{eq:final-value-condition}. Translating the reflected conclusion back gives $\Theta(b)\notin\PV$.
\end{proof}

\begin{claim}\label{clm:final-intervals}
We have
\(\ell-1<\Theta^{-1}(\ell)_y<\ell\) and \(u<\Theta(u)_x<u+1\).
\end{claim}

\begin{proof}
Because $h>\ell$ and $\Theta^{-1}(h)\in\PH$,
\(\Theta^{-1}(r)_y\le\Theta^{-1}(h)_y<\ell+1\).
Claim~\ref{clm:r-b-not-internal} excludes $\Theta^{-1}(r)$ from $\PH$, and $\PH$ is the only support set meeting the open strip $\ell<y<\ell+1$. Hence $\Theta^{-1}(r)_y\le \ell$, and therefore
\(\Theta^{-1}(\ell)_y<\Theta^{-1}(r)_y\le \ell\).
On the other hand, \eqref{eq:final-position-condition} places $\Sigma^{-1}(\ell-1)$ in $\Sigma_{\mathrm{BR}}$. Corollary~\ref{cor:BR-fixed} gives
\(\Theta^{-1}(\ell-1)=\Sigma^{-1}(\ell-1)\),
so $\Theta^{-1}(\ell-1)_y=\ell-1$ and
\(\Theta^{-1}(\ell)_y>\ell-1\).
This proves the first interval.

For the second, $v<u$ and \eqref{eq:h-v-ranges} give
\(\Theta(b)_x\ge\Theta(v)_x>u-1\), \(b=u-1\).
Claim~\ref{clm:r-b-not-internal} excludes $\Theta(b)$ from $\PV$, the only support set with $u-1<x<u$, so $\Theta(b)_x\ge u$. Hence
\(\Theta(u)_x>\Theta(b)_x\ge u\).
Finally, \eqref{eq:final-value-condition} places $\Sigma(u+1)$ in $\Sigma_{\mathrm{BR}}$, so Corollary~\ref{cor:BR-fixed} gives $\Theta(u+1)=\Sigma(u+1)$ and $\Theta(u+1)_x=u+1$. Thus $\Theta(u)_x<u+1$.
\end{proof}

We finish the proof. By Claim~\ref{clm:final-intervals},
\(\ell-1<\Theta^{-1}(\ell)_y<\ell\).
The only support sets meeting this horizontal unit strip are $\PR$ and $\PV$. The boundary restriction excludes $\Theta^{-1}(\ell)=\Theta(1)$ from $\PR$, so it lies in $\PV$. Since $\PV$ has $y$-coordinate range $(2,3)$ and $\ell\ge3$, this forces
\(\ell=3\), \(r=4\).

Similarly,
\(u<\Theta(u)_x<u+1\).
The only support sets meeting this vertical unit strip are $\PB$ and $\PH$. The boundary restriction excludes $\Theta(u)=\Theta^{-1}(k)$ from $\PB$, so it lies in $\PH$. Because $\PH$ has $x$-coordinate range $(k-2,k-1)$, we obtain
\(u=k-2\), \(b=k-3\).

Now \eqref{eq:final-value-condition} becomes
\(\sigma(2)<\sigma(u+1)<3\),
so $\sigma(2)=1$ and $\sigma(u+1)=2$. Since $b=\sigma^{-1}(1)=2$ and $b=k-3$, we get $k=5$ and $u=3$. The pattern is therefore
\(\sigma=(31524)\),
which was excluded because it was handled in Subsection~\ref{subsec:exceptional-five}. This final contradiction proves the lemma.
\end{proof}

Construction~\ref{con:two-internal-segments} covers the two stripes missing from Construction~\ref{con:four-boundary}, and Lemma~\ref{lem:two-internal-segments} gives uniqueness. Together with Corollaries~\ref{cor:large-gaps}, \ref{cor:horizontal-exceptions}, and~\ref{cor:vertical-exceptions}, and with the cases treated in Subsections~\ref{subsec:two-blocks}--\ref{subsec:coincident-boundary}, this completes the proof of Theorem~\ref{thm:geometric-construction} and hence of Theorem~1.2.

\paragraph*{Use of AI-assisted tools}

All mathematical results, constructions, and proofs in this paper were developed by the authors by the end of 2025. The results and proofs have already been presented and discussed in various seminars, and the authors had written a complete draft, before any AI tools were used for this manuscript. More recently during the last round of revision, the authors used ChatGPT and Claude to improve the exposition, to consolidate symmetric arguments, to help check notation, proof details, and references, and to assist with \LaTeX{} typesetting and with rendering figures designed and drawn by the authors. The authors independently reviewed and verified all revisions incorporated into the manuscript and take full responsibility for its content.

\paragraph{Formal verification.}
All numbered theorems, lemmas, claims, facts and corollaries of this paper, together with the
definitions and constructions we use, have been formalized in the Lean~4 proof
assistant~\cite{deMouraUllrich} using the Mathlib library~\cite{mathlib}, and our proofs checked by
Lean's kernel. The formalization contains no unproved steps and uses only Lean's standard axioms.
It is available at \url{https://github.com/ZRTMRH/le-shi-wei-linf-lean/tree/v1.0} (release v1.0,
which follows the numbering and statements of this version of the paper; Lean v4.35.0, Mathlib
commit \texttt{fcd79238}). Its file \texttt{CERTIFICATE.md} places the Lean statements of
Theorems~1.1, 1.2 and~4.1, and of the definitions we use, next to the statements in this paper;
the authors have reviewed and approved it. The Lean code was written by Claude using
\emph{formalize-paper}~\cite{formalizepaper}, an open-source Claude Code skill by Runtian Zhou,
which audits every statement with separate AI reviewers and produces a defect report and a
statement sign-off sheet for the authors.

\section*{Acknowledgment}
We thank Jacob Fox for a helpful discussion during the early stages of this project.
We thank Runtian Zhou for designing and running the workflow that produced the Lean formalization of this paper.


\end{document}